\documentclass{article}

\usepackage{amsthm}
\usepackage{amssymb}
\usepackage{amsmath}

\theoremstyle{definition}

\newtheorem{theorem}{Theorem}[section]
\newtheorem{proposition}[theorem]{Proposition}
\newtheorem{lemma}[theorem]{Lemma}
\newtheorem{corollary}[theorem]{Corollary}
\newtheorem{definition}[theorem]{Definition}

\theoremstyle{remark}
\newtheorem{remark}[theorem]{Remark}

\newcommand{\ihref}[1]{(\textbf{I\ref{#1}})}

\title{New Foundations is consistent}

\author{M. Randall Holmes and Sky Wilshaw}

\begin{document}

\maketitle

\begin{abstract}

In this paper we present a proof that Quine's set theory New Foundations of \cite{nf}, 1937 is consistent.
We do this by constructing a model of the ``tangled type theory" proposed by the first author and shown to be equiconsistent with NF in \cite{tangled}.  It is worth noting that the fact that the structure we describe is a model of tangled type theory has been formally verified in Lean by the second author and others (in \cite{wilshaw});  this paper alludes occasionally to the Lean formalization (and relevant sections are parallel to it in structure), but the paper version does not depend on the formal verification.

\end{abstract}

\tableofcontents

\newpage

\section{Introduction}

\subsection{Introductory note}

We are presenting an argument for the consistency of Quine's set theory New Foundations (NF) (see \cite{nf}).  The consistency of this theory relative to the usual systems of set theory has been an open question since the theory was proposed in 1937, with added urgency after Specker showed in 1954 that NF disproves the Axiom of Choice (see \cite{notac}).
Jensen showed in 1969 that NFU (New Foundations with extensionality weakened to allow urelements) is consistent and in fact consistent with Choice (and also with Infinity and with further strong axioms of infinity) (see \cite{nfu}).

The first author showed the equiconsistency of NF with the quite bizarre system TTT (Tangled Type Theory) in 1995 [see \cite{tangled}, but this paper is not really recommended; the presentation of its results here is better], which gave a possible approach to a consistency proof.  Since 2010, the first author has been attempting to write out proofs, first of the existence of a ``tangled web of cardinals" in Mac Lane set theory [not defined in this paper] and then directly of the existence of a model of tangled type theory.  These proofs have been  difficult to read, insanely involved, and involve the sort of elaborate bookkeeping which makes it easy to introduce errors every time a new draft is prepared.  We hope that the present version is more pleasant, though we know it is still forbidding.  The second author (with the assistance of others initially)
has formally verified the proof of the first author of the existence of a model of TTT (see \cite{wilshaw}), and so of the consistency of New Foundations, in the Lean proof verification system, and in the process has suggested minor corrections and considerable formal improvements to the argument originally proposed by the first author.  The second author reports that the formalized proof is still difficult to read and insanely involved with nasty bookkeeping.  Both authors feel that there ought to be a simpler approach, but the existing argument at least strongly resists attempts at simplification.

All the theories mentioned here are discussed (and referenced) in the next section.

Any remarks in the paper in the first person singular may be attributed to the first author.

\subsection{Acknowledgements}

The first author thanks Robert Solovay, who read a number of early versions of a related argument for Con(NF) and offered just criticisms which I have tried to take to heart.  I also thank Thomas Forster and Asaf Karagila, who have endured attempts of mine to present various relatives of this argument at length.  Further, I thank the members of the group of students at Cambridge who attempted formalization of an earlier version of this proof in Summer 2022.  We would like to thank Peter Lumsdaine for his work with the students in the summer project.   Cambridge students who participated in the summer project (including the second author) were funded partly by DPMMS, and Queens' College, but principally by the Wes and Margaret foundation, who had also supported work on New Foundations in the past, and we are grateful.

\newpage

\section{Development of relevant theories}

\subsection{The simple theory of types TST and TSTU}

We introduce a theory which we call the simple typed theory of sets or TST, a name favored by the school of Belgian logicians\footnote{led by Maurice Boffa, to whom the first author owes a profound debt, for being harbored in his group at the beginning of my career, which made it possible for me to launch my research program without distractions attendant on looking for a job.} who studied NF ({\em th\'eorie simple des types}).  This is not the same as the simple type theory of Ramsey and it is most certainly not Russell's type theory  (see historical remarks below).

TST is a first order multi-sorted theory with sorts (types) indexed by the nonnegative integers.  The primitive predicates of TST are equality and membership.

The type of a variable $x$ is written ${\tt type}(x)$:  this will be a nonnegative integer.   A countably infinite supply of variables of each type is supposed.   We provide a bijection $(x \mapsto x^+)$ from variables to variables of positive type satisfying   ${\tt type}(x^+)$ = ${\tt type}(x)+1$.\footnote{We do {\em not\/} furnish our variables with type superscripts.  One could follow the convention that if a variable has a natural number superscript, this determines its type, and that the effect of $^+$ on a superscripted variable is to increment the superscript, but we do not expect variables to have superscripts.}

An atomic equality sentence $x=y$ is well-formed iff ${\tt type}(x)={\tt type}(y)$. An atomic membership sentence $x \in y$ is well-formed iff ${\tt type}(x)+1 = {\tt type}(y)$.

The axioms of TST are extensionality axioms and comprehension axioms.

The extensionality axioms are all the well-formed assertions of the shape $(\forall xy:x=y \leftrightarrow (\forall z:z \in x \leftrightarrow z\in y))$.  For this to be well typed, the variables
$x$ and $y$ must be of the same type, one type higher than the type of $z$.

The comprehension axioms are all the well-formed assertions of the shape $(\exists A:(\forall x:x \in A \leftrightarrow \phi))$, where $\phi$ is any formula in which $A$ does not occur free.

The witness to $(\exists A:(\forall x:x \in A \leftrightarrow \phi))$ is unique by extensionality, and we introduce the notation $\{x:\phi\}$ for this object.  Of course, $\{x:\phi\}$  is to be assigned type one higher than that of $x$;  in general, complex terms will have types as variables do.

The modification which gives TSTU (the simple type theory of sets with urelements) replaces the extensionality axioms with the formulas of the shape $$(\forall xyw:w \in x \rightarrow (x=y \leftrightarrow (\forall z:z \in x \leftrightarrow z\in y))),$$  allowing many objects with no elements (called atoms or urelements)  in each positive type.  A technically useful refinement adds a constant $\emptyset^i$ of each positive type $i$ with no elements:  we can then address the problem that $\{x^i:\phi\}$ is not necessarily  uniquely defined when $\phi$ is uniformly false by defining $\{x^i:\phi\}$ as $\emptyset^{i+1}$ in this case.  The superscript is sometimes omitted from $\emptyset^i$ (and other constants parametrized by types) when the type can be deduced from context.

\subsubsection{Typical ambiguity}

TST(U) exhibits a symmetry, often called typical ambiguity, which is important in the sequel.

If $\phi$ is a formula, define $\phi^+$ as the result of replacing every variable $x$ (free and bound) in $\phi$ with $x^+$ (and occurrences of $\emptyset^i$ with $\emptyset^{i+1}$ if this is in use).   It should be evident that if $\phi$ is well-formed, so is $\phi^+$,
and that if $\phi$ is a theorem, so is $\phi^+$ (the converse is not the case).  Further, if we define a mathematical object as a set abstract $\{x:\phi\}$ we have an analogous
object $\{x^+:\phi^+\}$ of the next higher type (this process can be iterated).

The axiom scheme asserting $\phi \leftrightarrow \phi^+$ for each closed formula $\phi$ is called the Ambiguity Scheme.   Notice that this is a stronger assertion than is warranted by the symmetry of proofs described above.

\subsubsection{Historical remarks}

TST is not the type theory of the {\em Principia Mathematica\/} of Russell and Whitehead (\cite{pm}), though a description of TST is a common careless description of Russell's theory of types.

Russell described something like TST informally in his 1904 {\em Principles of Mathematics\/} (\cite{pm1}).  The obstruction to giving such an account in {\em Principia Mathematica\/} was that
Russell and Whitehead did not know how to describe ordered pairs as sets.  As a result, the system of {\em Principia Mathematica\/} has an elaborate system of  complex
types inhabited by $n$-ary relations with arguments of specified previously defined types, further complicated by predicativity restrictions (which are in effect cancelled by an axiom of reducibility).
The simple theory of types of Ramsey eliminates the predicativity restrictions and the axiom of reducibility, but is still a theory with complex types inhabited by $n$-ary relations.

Russell noticed a phenomenon like the typical ambiguity of TST in the more complex system of {\em Principia Mathematica\/}, which he refers to as ``systematic ambiguity".

In 1914 (\cite{wiener}), Norbert Wiener gave a definition of the ordered pair as a set (not the one now in use) and seems to have recognized that the type theory of {\em Principia Mathematica\/} could be simplified to something like TST, but he did not give a formal description.  The theory we call TST was apparently first described by Tarski in the 1930s\footnote{We regard it as an important feature of TST that the nature of type 0 is left completely undetermined, so we do not regard descriptions of arithmetic of order $\omega$ appearing about the same time as part of this history.     Arithmetic of order $\omega$ is a similar typed theory but lacks the ambiguity of interest here.}.

It is worth observing that the axioms of TST look exactly like those of ``naive set theory", the restriction preventing paradox being embodied in the restriction of the language by the type system.
For example, the Russell paradox is averted because one cannot have $\{x:x \not\in x\}$ because $x \in x$ (and so its negation $\neg x \in x$) cannot be a well-formed formula.

It was shown around 1950 (in \cite{kemeny}) that Zermelo set theory proves the consistency of TST with the axiom of infinity;  TST + Infinity has the same consistency strength as
Zermelo set theory with separation restricted to bounded formulas.

\newpage

\subsubsection{Some mathematics in TST}

We briefly discuss some mathematics in TST.

We indicate how to define the natural numbers.  We use the definition of Frege ($n$ is the set of all sets with $n$ elements).  0 is $\{\emptyset\}$ (notice that we get a natural number 0 in each type $i+2$;  we will be deliberately ambiguous in this discussion, but we are aware that anything we define is actually not unique, but reduplicated in each type above the lowest one in which it can be defined).  For any set $A$ at all we define $\sigma(A)$ as $\{a \cup \{x\}:a \in A \wedge x \not\in a\}$.  This is definable for any $A$ of type $i+2$ ($a$ being of type $i+1$ and $x$ of type $i$).  Define 1 as $\sigma(0)$, 2 as $\sigma(1)$,  3 as $\sigma(2)$, and so forth.  Clearly we have successfully defined 3 as the set of all sets with three elements, without circularity.
But further, we can define $\mathbb N$ as $\{n:(\forall I:0 \in I \wedge (\forall x \in I:\sigma(x) \in I) \rightarrow n \in I)\}$, that is, as the intersection of all inductive sets.
$\mathbb N$ is again a typically ambiguous notation:  there is an object defined in this way in each type $i+3$.

The collection of all finite sets can be defined as $\bigcup \mathbb N$.  The axiom of infinity can be stated as $V \not\in \bigcup \mathbb N$ (where $V= \{x:x=x\}$ is the typically ambiguous symbol for the type $i+1$ set of all type $i$ objects).  It is straightforward to show that the natural numbers in each type of a model of TST with Infinity are isomorphic in a way representable in the theory.

Ordered pairs can be defined following Kuratowski and a quite standard theory of functions and relations can be developed.  Cardinal and ordinal numbers can be defined as Frege or Russell would have defined them, as isomorphism classes of sets under equinumerosity and isomorphism classes of well-orderings under similarity.

The Kuratowski pair $(x,y) = \{\{x\},\{x,y\}\}$ is of course two types higher than its projections, which must be of the same type.  There is an alternative definition (due to Quine in \cite{quinepair}) of an ordered pair
$\left< x,y\right>$ in TST + Infinity which is of the same type as its projections $x,y$.  This is a considerable technical convenience but we will not need to define it here.  Note for example that if we use the Kuratowski pair, the cartesian product $A \times B$ is two types higher than $A,B$, so we cannot define $|A| \cdot |B|$ as $|A \times B|$ if we want multiplication of cardinals to be a sensible operation.  Let $\iota$ be the singleton operation and define $T(|A|)$ as $|\iota``A|$ (this is a very useful operation sending cardinals of a given type to cardinals in the next higher type which seem intuitively to be the same; also, it is clearly injective, so has a (partial) inverse operation $T^{-1}$).  The definition of cardinal multiplication if we use the Kuratowski pair is then $|A| \cdot |B| =T^{-2}(|A\times B|)$.  If we use the Quine pair this becomes the usual definition $|A| \cdot |B| =|A\times B|$.  Use of the Quine pair simplifies matters in this case, but it should be noted that the $T$ operation remains quite important (for example it provides the internally representable isomorphism between the systems of natural numbers in each sufficiently high type).

Note that the form of Cantor's Theorem in TST is not $|A| < |{\cal P}(A)|$, which would be ill-typed, but $|\iota``A|<|{\cal P}(A)|$:  a set has fewer unit subsets than subsets.  The exponential map $\exp(|A|) = 2^{|A|}$ is not defined as $|{\cal P}(A)|$, which would be one type too high, but as $T^{-1}(|{\cal P}(A)|)$, the cardinality of a set $X$ such that $|\iota``X| = |{\cal P}(A)|$;   notice that this is partial.  For example
$2^{|V|}$ is not defined (where $V=\{x:x=x\}$, an entire type), because there is no $X$ with $|\iota``X|=|{\cal P}(V)|$, because $|\iota``V|<|{\cal P}(V)| \leq |V|$, and of course there is no set larger than $V$ in its type.

\newpage

\subsection{New Foundations and NFU}

In \cite{nf}, 1937, Willard van Orman Quine proposed a set theory motivated by the typical ambiguity of TST described above.  The paper in which he did this was titled ``New foundations for mathematical logic", and the set theory it introduces is called ``New Foundations" or NF, after the title of the paper.\footnote{There is a persistent rumor that Quine's set theory in its original version was inconsistent.  This is not the case (and in fact we show here that the system of the 1937 paper in which the system was introduced is consistent).  The truth behind the rumor is that the system of the first edition of Quine's book {\em Mathematical Logic\/} (\cite{ml}), in which proper classes were added to NF, was inconsistent;  the inconsistency was corrected in the second edition of 1951, and the results of this paper show that the system of the 1951 book is consistent.  We do not discuss the systems with proper classes further here.}

Quine's observation is that since any theorem $\phi$ of TST is accompanied by theorems $\phi^+, \phi^{++}, \phi^{+++}, \ldots$ and every defined object $\{x:\phi\}$ is accompanied by
$\{x^+:\phi^+\},\{x^{++}:\phi^{++}\},\{x^{+++}:\phi^{+++}\}$, so the picture of what we can prove and construct in TST looks rather like a hall of mirrors\footnote{I have been asked whether the turn of phrase ``hall of mirrors" here is my own or due to Quine himself or others in the interim:  I actually do not know!}, we might reasonably (?) suppose that the types are all the same.

The concrete implementation follows.  NF is the first order unsorted theory with equality and membership as primitive with an axiom of extensionality $(\forall xy:x=y \leftrightarrow (\forall z:z \in x \leftrightarrow z\in y))$ and an axiom of comprehension $(\exists A:(\forall x:x \in A \leftrightarrow \phi))$ for each formula $\phi$ in which $A$ is not free which can be obtained from a formula of TST by dropping all distinctions of type.  We give a precise formalization of this idea:  provide a bijective map $(x \mapsto x^*)$ from the countable supply of variables (of all types) of TST onto the countable supply of variables of the language of NF.  Where $\phi$ is a formula of the language of TST, let $\phi^*$ be the formula obtained by replacing every variable $x$, free and bound,
in $\phi$ with $x^*$. For each formula $\phi$ of the language of TST in which $A$ is not free in $\phi^*$ and each variable $x^*$, an axiom of comprehension of NF asserts $(\exists A:(\forall x^*:x^* \in A \leftrightarrow \phi^*))$.

In the original paper, this is expressed in a way which avoids explicit dependence on the language of another theory.  We do this in a way which is usual now but not necessarily identical to what is done in the original paper.  Let $\phi$ be a formula of the language of
NF.    A function $\sigma$ is a stratification of $\phi$ if it is a (possibly partial) map from variables to non-negative integers such that for each atomic subformula
`$x=y$'  of $\phi$ we have $\sigma($`$x$'$)=\sigma($`$y$'$)$ and for each atomic subformula `$x \in y$' of $\phi$ we have $\sigma($`$x$'$)+1 = \sigma($`$y$'$)$.
A formula $\phi$ is said to be stratified iff there is a stratification of $\phi$.  Then for each stratified formula $\phi$ of the language of NF and variable $x$ we have an axiom $(\exists A:(\forall x:x \in A \leftrightarrow \phi))$.  The stratified formulas are exactly the formulas $\phi^*$ up to renaming of variables.

For each stratified formula $\phi$, there is a unique witness to $$(\exists A:(\forall x:x \in A \leftrightarrow \phi))$$ (uniqueness follows by extensionality) which we denote by $\{x:\phi\}$.

NF has been dismissed as a ``syntactical trick" because of the way it is defined.  It might go some way toward dispelling this impression to note that the stratified comprehension scheme is equivalent to a finite collection of its instances, so the theory can be presented in a way which makes no reference to types at all.  This is a result of Hailperin (\cite{hailperin}), refined by others, and discussed in some detail in a subsequent subsection of this section, because it is relevant to the Lean formalization.  One obtains a finite axiomatization of NF by analogy with the method of finitely axiomatizing von Neumann--G\"odel--Bernays predicative class theory.  It should further be noted that the first thing one does with any finite axiomatization is prove stratified comprehension as a meta-theorem, in practice, but it remains significant that the theory can be axiomatized with no reference to types at all.  It is also worth noting that the collection of all typed instances of the Hailperin axioms is an axiomatization of TST, not a finite axiomatization, but an axiomatization with finitely many templates.

Jensen in \cite{nfu}, 1969 proposed the theory NFU which replaces the extensionality axiom of NF with $$(\forall xyw:w \in x \rightarrow (x=y \leftrightarrow (\forall z:z \in x \leftrightarrow z\in y))),$$  allowing many atoms or urelements.  One can reasonably add an elementless constant $\emptyset$, and define $\{x:\phi\}$ as $\emptyset$ when $\phi$ is false for all $x$.

Jensen showed that NFU is consistent and moreover NFU + Infinity + Choice is consistent.  We will give an argument similar in spirit though not the same in detail for the consistency of NFU in the next section.

An important theorem of Specker (\cite{ambiguity}, 1962) is that NF is consistent if and only if TST + the Ambiguity Scheme is consistent.  His method of proof adapts to show that  NFU is consistent if and only if TSTU + the Ambiguity Scheme is consistent.  Jensen used this fact in his proof of the consistency of NFU.  We prove a version of Specker's result using concepts from this paper below.

In \cite{notac}, 1954, Specker had shown that NF disproves Choice, and so proves Infinity.  At this point if not before it was clear that there is a serious issue of showing that NF is consistent relative to some set theory in which we have confidence.  There is no evidence that NF is any stronger than TST + Infinity, the lower bound established by Specker's result.

Note that NF or NFU supports the implementation of mathematics in the same style as TST, but with the representations of mathematical concepts losing their ambiguous character.  The number 3 really is realized as the unique set of all sets with three elements, for example.  The universe is a set and sets make up a Boolean algebra.   Cardinal and ordinal numbers can be defined
in the manner of Russell and Whitehead.

The apparent vulnerability to the paradox of Cantor is an illusion.  Applying Cantor's theorem to the cardinality of the universe in NFU gives $|\iota``V| < |{\cal P}(V)| \leq |V|$ (the last inequality would be an equation in NF), from which we conclude that there are fewer singletons of objects than objects in the universe.  The operation $(x \mapsto \{x\})$ is not a set function, and there is every reason to expect it not to be, as its definition is unstratified.  The resolution of the Burali-Forti paradox is also weird and wonderful in NF(U), but would take us too far afield.

\newpage

\subsection{Tangled type theory TTT and TTTU}

In \cite{tangled}, 1995, the first author described a reduction of the NF consistency problem to consistency of a typed theory,  motivated by reverse engineering from Jensen's method of proving the consistency of NFU.

Let $\lambda$ be a limit ordinal.  It can be $\omega$ but it does not have to be.

In the theory TTT (tangled type theory) which we develop, each variable $x$ is supplied with a type ${\tt type}($`$x$'$) <\lambda$;  we are provided with countably many distinct variables of each type.

For any formula $\phi$ of the language of TST and any strictly increasing sequence $\{s_i\}_{i \in \mathbb N}$ in $\lambda$, let $\phi^s$ be the formula obtained by replacing each variable
of type $i$ with a variable of type $s(i)$.  To make this work rigorously, we suppose that we have a bijection from type $i$ variables of the language of TST to type $\alpha$ variables
of the language of TTT for each natural number $i$ and ordinal $\alpha<\lambda$.

TTT is then the first order theory with types indexed by the ordinals below $\lambda$ whose well-formed atomic sentences `$x=y$' have ${\tt type}($`$x$'$) = {\tt type}($`$y$'$)$ and whose well-formed atomic sentences `$x \in y$' satisfy ${\tt type}($`$x$'$) < {\tt type}($`$y$'$)$, and whose axioms are the sentences $\phi^s$ for each axiom $\phi$ of TST and each strictly increasing sequence $s$ in $\lambda$.  TTTU has the same relation to TSTU (with the addition of constants $\emptyset^{\alpha,\beta}$ for each $\alpha<\beta<\lambda$  such that $(\forall x^{\alpha} :x^{\alpha}\not\in \emptyset^{\alpha,\beta})$ is an axiom).

It is important to notice how weird a theory TTT is.  This is not cumulative type theory.  Each type $\beta$ is being interpreted as a power set of {\em each\/} lower type $\alpha$.  Cantor's theorem in the metatheory makes it clear that most of these power set interpretations cannot be honest.

There is now a striking

\begin{theorem}[Holmes]
  TTT(U) is consistent iff NF(U) is consistent.
\end{theorem}
\begin{proof}
Suppose NF(U) is consistent.  Let $(M,E)$ be a model of NF(U) (a set $M$ with a membership relation $E$).  Implement type $\alpha$ as $M \times \{\alpha\}$ for
each $\alpha<\lambda$.  Define $E_{\alpha,\beta}$ for $\alpha<\beta$ as $\{((x,\alpha),(y,\beta)):xEy\}$.  This gives a model of TTT(U).   Empty sets in TTTU present no essential additional difficulties.

Suppose TTT(U) is consistent, and so we can assume we are working with a fixed model of TTT(U).  Let $\Sigma$ be a finite set of sentences in the language of TST(U).  Our first aim is to construct a model of TST(U) in which $\phi \leftrightarrow \phi^+$ holds for each sentence $\phi \in \Sigma$.
\begin{enumerate}
\item Let $n$ be the smallest natural number such that all variables which occur in $\Sigma$ have type less than $n$. We define a partition of the $n$-element subsets of $\lambda$.  Each $A \in [\lambda]^n$ is put in a compartment
determined by the truth values of the sentences $\phi^s$ in our model of TTT(U), where $\phi \in \Sigma$ and ${\tt rng}(s \restriction \{0,\ldots,n-1\}) = A$:  there are no more than $2^{|\Sigma|}$ compartments.  By Ramsey's theorem, there is an infinite set $H \subseteq \lambda$ homogeneous for this partition  which includes the range of a strictly increasing sequence $h$.

\item Any strictly increasing sequence $s$ determines a model of TST(U) in which each type $i$ is implemented as type $s(i)$ of our fixed model of TTT(U)  and membership in the derived model of TST(U) of type $i$ objects in type $i+1$ objects is implemented as membership in the fixed model of TTT(U) of type $s(i)$ objects in type $s(i+1)$ objects.

\item This model determined in this way by $h$ satisfies  $\phi \leftrightarrow \phi^+$ for each $\phi \in \Sigma$.   This achieves our first aim.

\item But this implies by compactness that the full Ambiguity Scheme $\phi \leftrightarrow \phi^+$ is consistent with TST(U).  We could at this point say that  NF(U) is consistent by the 1962 result of Specker (\cite{ambiguity}), but we will continue and indicate how to prove this directly (our second aim).

\item We pause to add a useful notion to our logic.  The Hilbert symbol is a primitive term construction $(\epsilon x:\phi)$ (same type as $x$) with axiom schemes $(\exists x:\phi) \rightarrow \phi[(\epsilon x:\phi)/x]$ and $(\forall x:\phi \leftrightarrow \psi) \rightarrow (\epsilon x:\phi) = (\epsilon x:\psi)$ which cannot appear in instances of comprehension. We do not define the quantifiers in terms of the Hilbert symbol, because the quantifiers do need to appear in instances of comprehension.

\item If we have a model of TST(U) augmented with a Hilbert symbol  which satisfies Ambiguity (for all formulas, including those which mention the Hilbert symbol) then we can readily get a model of NF(U), by constructing a term model using the Hilbert symbol in the natural way, then identifying all terms with their type-raised versions (an element of the term model is a class of terms $(\epsilon x:\phi)$ in the language of TST(U) which is obtained from a single such term by all possible upward or downward shifts of type by a constant amount:  we write this
$[(\epsilon x:\phi)]$).

\item All statements in the resulting type-free theory can be decided:  the truth value of a closed atomic formula $[(\epsilon x:\phi)] \,R\, [(\epsilon y:\psi)]$ in the model of NF(U) is determined by replacing each side of the formula with one of its elements in such a way as to get a well-typed formula, which will always be possible if types are taken high enough, and reading the truth value of the resulting formula  from the model of TST(U) with Ambiguity and Hilbert symbol;  $R$ is either = or $\in$ [note that the relation interpreting equality in the term model is a nontrivial equivalence relation].  Truth values of more complex statements are determined in standard ways.  That this makes axioms of NF(U) true should be clear.

\item  Now observe that a model of TTT(U) can readily be equipped with a Hilbert symbol if this creates no obligation to add instances of comprehension
containing the Hilbert symbol (use a well-ordering of the set implementing each type to interpret a Hilbert symbol  $(\epsilon x:\phi)$ in that type as the first $x$ such that $\phi$ or as a default object of the appropriate type if there is no such $x$), and the argument in paragraphs (1)--(4) adapts to show consistency of TST(U) plus Ambiguity with the Hilbert symbol.

\item  We remark that our method of proving (a version of) the Specker ambiguity result is distinctly different from the original method.  Strictly speaking, we have only proved Specker's result for TST(U) + Ambiguity + existence of a Hilbert symbol, but that is all we need for our result.

\item From this it follows that NF(U) is consistent under the hypothesis that TTT(U) is consistent, which is what was to be shown.

\end{enumerate}
\end{proof}

\begin{theorem}[essentially due to Jensen]
  NFU is consistent.
\end{theorem}
\begin{proof}
It is enough to exhibit a model of TTTU.  Represent type $\alpha$ as $V_{\omega+\alpha} \times \{\alpha\}$ for each $\alpha<\lambda$ ($V_{\omega+\alpha}$ being the rank indexed by $\omega+\alpha$ in the usual cumulative hierarchy).  Define $\in_{\alpha,\beta}$ for
$\alpha<\beta<\lambda$ as $$\{((x,\alpha),(y,\beta)):x \in V_{\omega+\alpha} \wedge y \in V_{\omega+\alpha+1} \wedge x \in y\}.$$  This gives a model of TTTU in which the membership of
type $\alpha$ in type $\beta$ is defined in such a way that each $(y,\beta)$ with $y \in V_{\omega+\beta} \setminus V_{\omega+\alpha+1}$ is interpreted as an urelement over type $\alpha$.

Our use of $V_{\omega+\alpha}$ enforces Infinity in the resulting models of NFU (note that we did not have to do this:  if we set $\lambda=\omega$ and interpret type $\alpha$ using $V_\alpha$ we prove the consistency of NFU with the negation of Infinity).  It should be clear that Choice (stated in a stratified way) holds in the models of NFU eventually obtained if it holds in the ambient set theory.

This shows in fact that mathematics in NFU is quite ordinary (with respect to stratified sentences), because mathematics in the models of TSTU embedded in the indicated model of TTTU is quite ordinary.  The notorious ways in which NF evades the paradoxes of Russell, Cantor and Burali-Forti can be examined in actual models and we can see that they work and how they work (since they work in NFU in the same way they work in NF).
\end{proof}

Of course Jensen did not phrase his argument in terms of tangled type theory.  Our contribution here was to reverse engineer from Jensen's original argument for the consistency of NFU an argument for the consistency of NF itself, which requires additional input which we did not know how to supply (a proof of the consistency of TTT itself).  An intuitive way to say what is happening here is that Jensen noticed that it is possible to skip types in a certain sense in TSTU in a way which is not obviously possible in TST itself;  to suppose that TTT might be consistent is to suppose that such type skipping is also possible in TST, and the description of TTT is an explicit realization of what is meant by this.  We leave it to the historically inclined reader to examine Jensen's original proof and determine the exact relationship between the approaches.

\newpage

\subsubsection{There are $\alpha$-models for every $\alpha$}\label{subsection:alpha_models}

This section is devoted to adapting a result of Jensen to our context.  Jensen proved in \cite{nfu} that for any infinite ordinal $\alpha$ of the metatheory, there is a model of NFU + Infinity in which there is a well-ordering of type $\alpha$.  This is of particular interest in the case $\alpha=\omega$, which establishes that the natural numbers of a model can correspond precisely to the natural numbers of the metatheory.

Nothing in our construction of a model of TTT below depends on this section:  the fact that an existence of an $\omega$-model of NF follows from our main result is significant for extensions of our results and uses this theorem.

A relation $R$ of one of the theories under discussion is a set of ordered pairs\footnote{The discussion here presupposes use of the Kuratowski ordered pair;  it could be adapted to the use of another pair definition.}  of the theory in question;  note that it corresponds to a set of ordered pairs of the metatheory with the same field.  If this relation corresponding to  $R$ is a well-ordering, we can talk about its order type both in terms of the theory under consideration and in terms of the metatheory, though ordinals of the theory under consideration and the metatheory may not be the same.

In any of the theories under discussion, an ordinal is taken to be an equivalence class of well-orderings under similarity.
The relations which belong to an ordinal of one of these theories may fail to be well-orderings in terms of the metatheory.
If $\alpha$ is an ordinal of the metatheory, we say that $\alpha$ is {\em implemented\/} in a type of one of the theories under discussion
exactly if there is an ordinal in that type (in the sense of the theory) whose elements in fact correspond to well-orderings of order type $\alpha$ in the sense of the metatheory, in the sense outlined in the previous paragraph.  It is straightforward to show that if $\alpha$ is implemented in a type, the ordinal implementing it will contain all well-orderings of appropriate type in the sense of the theory  whose order type is $\alpha$ in the sense of the metatheory.

\begin{theorem}

For any infinite ordinal $\alpha$ of the metatheory, if there is a model of TTT(U) in which the ordinal $\alpha$ has an implementation $\alpha^s$ for each strictly increasing sequence $s$ of types, by which we mean an ordinal in type $s_4$ whose elements of type $s_3$ correspond to well-orderings of order type $\alpha$ of subsets of type $s_0$, using types $s_1$ and $s_2$ for the intermediate pairing machinery, and which is such that the ordinal $\lambda$ indexing the types is a strong limit cardinal of cofinality $>|\alpha|^{|\alpha|}$, then there is a model of NF(U) in which the ordinal $\alpha$ is implemented, which further satisfies [after neglect of types] any sentence in the language of TST(U) whose translations via all increasing sequences of types are satisfied in the model of TTT(U).

\end{theorem}

\begin{proof}

The proof is a modified version of the proof given above for consistency of NF(U), using a more powerful partition principle.  We do not claim originality, referring the reader to \cite{nfu}, but our presentation is not the same, and it is adapted to talk of tangled type theory.

The Erd\H{o}s--Rado theorem tells us that for a natural number $n\geq 1$ and an infinite cardinal $\mu$, any partition of $\exp^n(\mu)$ into $\mu$ compartments
has a homogeneous set of size $\mu^+$.  This is not quite the usual formulation, but easily follows from it.

We extend the language of TTT(U) to include constant terms $\alpha^s$ and $\beta^s$ for each ordinal $\beta < \alpha$ of the metatheory for each sequence $s$ as described above.   We similarly extend the language of TST(U) (in this case only needing constants $\alpha^i$ and $\beta^i$ for $\beta<\alpha$ for each $i\geq4$),  and note that the translation of a sentence of TST(U) with language augmented with these ordinal constants into terms of a strictly increasing sequence of types in TTT(U) is straightforward, and we allow use of these constants in the comprehension axiom of TTT(U) as usual.  We further include the Hilbert symbol in our language, with the same restriction stated above that the Hilbert symbol may not appear in instances of comprehension, noting our intention to construct a term model.

We consider partitions of $[\lambda]^n$  determined by  the set $\Sigma_n$ of {\em all\/} closed Hilbert symbols (``closed" meaning the symbol contains no parameters) of the form $(\epsilon \beta^i:\beta^i < \alpha^i \wedge \phi)$
in the language of TST(U) augmented  as indicated above and using only types less than $n$, with the compartment in the partition containing $A \in [\lambda]^n$ determined by the ordinal values $<\alpha$ [or default value, where we assume the default value is not an ordinal] assigned to all of these terms in our model of TTT(U) when types $0-(n-1)$ are replaced by the types in $A$ in increasing order in each of the Hilbert symbols in $\Sigma_n$.  Here we are speaking loosely in that we say an ordinal $\gamma<\alpha$ is assigned as value to $(\epsilon \beta^i:\beta^i < \alpha^i \wedge \phi)$ when its value is in fact $\gamma^i$:  it should be clear that this is unproblematic and easier to say.

There are no more than $|\alpha|$ such Hilbert symbols, and no more than $|\alpha|^{|\alpha|}$ compartments in each of these partitions.

Note that assignment of values to all such Hilbert symbols induces a complete assignment of truth values to sentences of the extended language:  $$(\epsilon \beta^i:\beta^i<\alpha \wedge \beta^i<2 \wedge (\beta=1 \leftrightarrow \phi))$$ codes a sentence $\phi$ into $\Sigma_n$ if $\phi$ mentions only types below $n > i$, because the assignment of ordinal value to this term corresponds to the assignment of truth value to $\phi$.

The Erd\H{o}s--Rado theorem allows us to conclude that for any infinite $\mu<\lambda$, we have a homogeneous set of size $\mu^+$ for the given partition in $\exp^n(\mu)<\lambda$.
Further, we can find a single assignment of values to Hilbert symbols as described above which is associated with arbitrarily large homogeneous sets below $\lambda$, because the cofinality of $\lambda$ is greater than $|\alpha|^{|\alpha|}$, the cardinality of the set of compartments in the partition, corresponding precisely to the cardinality of the set of such assignments of values.
If we have chosen such an assignment of values for a specific value of $n$, we can further choose such an assignment of values for $n+1$ which extends it:  for any $\mu$, choose a homogeneous set for the partition of $[\lambda]^n$ realizing the fixed assignment of values and of cardinality $\nu > \exp^{n+1}(\mu)$, then find a homogeneous set of size $\mu^+$ for the partition of $[\lambda]^{n+1}$ inside this homogeneous set, which will determine an assignment of values to Hilbert symbols using $n+1$ types extending the fixed one chosen for Hilbert symbols using $n$ types.  Thus there are arbitrarily large homogeneous sets for
assignments of values to Hilbert symbols using $n+1$ types extending the fixed one for Hilbert symbols using $n$ types.  Thus there is a fixed assignment of values to Hilbert symbols using $n+1$ types which extends the given fixed assignment for $n$ types and which is associated with arbitrarily large homogeneous sets, by the already stated considerations of cofinality.  Thus we can proceed through $\omega$ steps and find an assignment of values to all Hilbert symbols with the property that its restriction to any fixed number $n$ of types is realized in arbitrarily large homogenous sets below $\lambda$.

Thus we have constructed a complete assignment of values to such Hilbert symbols in the extended TST(U) language, and thus a complete assignment of truth values to all sentences, which is ambiguous [invariant under displacement of types] and consistent with any statement holding at all sequences of types in the model of TTT(U).  Thus for example we will get extensionality if the model of TTTU is actually a model of TTT.   In fact, for any fixed $n$ all these assignments for Hilbert symbols mentioning only types less than $n$ will be realized at some sequence of types in the model of TTT(U) [verifying that this assignment of values is consistent and in effect a complete description of a model of NF(U)].

Further, this model of NFU contains no ordinal less than what it calls $\alpha$ other than the ordinals which it calls $\beta$ for $\beta<\alpha$ in the metatheory, and the natural well-ordering on these ordinals is realized in the model of NF(U) and corresponds to a well-ordering of order type $\alpha$ in the metatheory.  The reason for this is that the model of NFU is a term model, and any term denoting an ordinal less than $\alpha$ can be slightly modified to the exact form used in defining the partitions, and the process will assign that term a concrete value less than $\alpha$, an actual ordinal $\beta<\alpha$ of the metatheory, so the order on the ordinals $<\alpha$ in the model of NF(U) has order type $\alpha$ in the metatheory.

\end{proof}

\newpage

\subsection{An axiomatization of TST with finitely many templates}\label{ss:hailperin}

We discuss a finite axiomatization of NF derived from that of Hailperin (the statement of the axioms is taken from an implementation of Hailperin's axiom set of \cite{hailperin} in Metamath (\cite{metamath}), and there are minor changes from the original formulation;  the proofs are our own), making the important observation that it actually provides us with an axiomatization of TST with finitely many
axiom templates (in the sense that each axiom is a type shifted version of one of a finite set of axioms).  Notation introduced in this section is not used in the rest of the paper, and nothing in subsequent sections except a brief remark in the last paragraphs of section 4 depends on anything here.  This finite axiomatization is however used in the Lean formalization, which is why we review it in detail here.

The finite axiomatization of NF takes this form (the definitions inserted are ours, and we have modified the order of the axioms to make the definitions work sensibly).  We also present this as an axiomatization of TST with finitely many templates, with the proviso that each typed form of each axiom is asserted:

\begin{description}

\item[extensionality axiom:]  $(\forall x:(\forall y:(\forall z:z \in x \leftrightarrow z \in y)\rightarrow x=y))$

\item[anti-intersection axiom:]  $(\forall xy:(\exists z:(\forall w:w \in z \leftrightarrow \neg(w \in x \wedge w \in y))))$

\item[singleton axiom:]  $(\forall x:(\exists y:(\forall z:z \in y \leftrightarrow z = x)))$

\item[definition:]  $\{x\}$ denotes for each $x$ the set whose only element is $x$, whose existence is provided by the previous axiom.  We define $\iota(x)$ as $\{x\}$ and define $\iota^1(x)$ as $\iota(x)$ and $\iota^{n+1}(x)$ as $\{\iota^n(x)\}$, for each concrete natural number $n$.

\item[cardinal one axiom:] $(\exists x:(\forall y:y \in x \leftrightarrow (\exists z:(\forall w:w \in y \leftrightarrow w = z))))$

\item[definition:]  We define 1 as the set of all singletons, provided by the previous axiom.

\item[definition:]   $x | y$ denotes the set $z$ whose existence is provided by the anti-intersection axiom:  $z \in x | y \leftrightarrow \neg(z \in x \wedge z \in y)$.
We define $x^c$ as $x | x$.  We define $x \cap y$ as $(x|y)^c$.  We define $x \cup y$ as $x^c | y^c$.  We define $V$ as $1|1^c$ (noting that it is straightforward to prove $x|x^c = V$ for any $x$, since this is the universal set).  We define $\{x,y\}$ as $\{x\} \cup \{y\}$.  We define $(x,y)$ as $\{\{x\},\{x,y\}\}$.  We define
$(x,y,z)$ as $(\{\{x\}\},(y,z))$.   More generally, we define $(x_1,\ldots,x_n)$ as $(\iota^{2n-4}(x_1),(x_2,\ldots,x_n))$ for $n>2$ a concrete natural number. [The treatment of $n$-tuples is what makes this axiomatization singularly awful].

\item[cross product axiom:]  $(\forall x:(\exists y:(\forall z:z \in y \leftrightarrow (\exists wt:z = (w,t) \wedge t \in x))))$

\item[definition:]  We define $V \times x$ as the set introduced by the previous axiom:  $z \in V  \times x \leftrightarrow (\exists wt:z = (w,t) \wedge t \in x)$.  Note that $V \times V$ is the set of all ordered pairs.

\item[converse axiom:] $(\forall x:(\exists y:(\forall zw:(z,w)\in y \leftrightarrow (w,z) \in x)))$

\item[definition:]  For any set $R$, we define $R^{-1}$ as the intersection of $V \times V$ with a set introduced by the previous axiom:  $$(\forall zw:(z,w)\in R^{-1} \leftrightarrow (w,z) \in R) \wedge (\forall u:u \in R^{-1} \rightarrow (\exists zw:(z,w)=u)).$$

\item[definition:]  We define $x \times V$ as $(V \times x)^{-1}$ and $x \times y$ as $(x \times V) \cap (V \times y)$.

\item[definitions:] We define $\iota^2``x$ as $(x \times x) \cap 1$.  This is the image of $x$ under the double singleton operation.

Note that $\iota^2``V$ is the equality relation.  Define $\iota^{2(n+1)}``x$ as $\iota^2``(\iota^{2n}``x)$.

We define $x_1 \times x_2 \ldots \times x_n$ as $\iota^{2(n-2)}``x \times (x_2 \times \ldots \times x_n)$.

We define $V^1$ as $V$ and $V^{n+1}$ as $\iota^{2(n-1)}``V\times V^n$, for each concrete $n$.  $V^n$ is the set of all $n$-tuples.

\item[singleton image axiom:]  $(\forall x:(\exists y:(\forall zw:(\{z\},\{w\}) \in y \leftrightarrow (z,w) \in x)))$.
\item[definition:]  We define $R^{\iota}$ for any set $R$ as the intersection of a set provided by the previous axiom with $1 \times 1$.  $R^{\iota}$ is the set whose members are exactly the ordered pairs $(\{z\},\{w\})$ such that $(z,w)\in R$.  Let $R^{\iota^1}$ be defined as $R^{\iota}$ and $R^{\iota^{n+1}}$ be defined as $(R^{\iota^n})^\iota$.

\item[insertion two axiom:]  $(\forall x:(\exists y:(\forall zwt:(z,w,t) \in y \leftrightarrow (z,t) \in x)))$

We define $I_2(x)$ as the intersection of a set provided by the previous axiom with $V \times V \times V$.

\item[insertion three axiom:] $(\forall x:(\exists y:(\forall zwt:(z,w,t) \in y \leftrightarrow (z,w) \in x)))$

We define $I_3(x)$ as the intersection of a set provided by the previous axiom with $V \times V \times V$.

\item[definition:]  It seems natural to define $I_1(R)$ as $\iota^2``V \times R$; this requires no new axiom.

\item[definition:]  Define $I_{1,n}(R)$ as $\iota^{2(n-1)}``V \times R$:  this is correct for prepending all possible initial projections
to an $n$-tuple.  Then define $I^1_{1,n}(R)$ as $I_{1,n}(R)$ and define $I^{m+1}_{1,n}(R)$ as $I_{1,n+m}(I^m_{1,n}(R))$:  this takes
into account the fact that the tuples get longer.   The intention here is that $I^m_{1,n}$ represents the result of applying the insertion one operation to an $n$-ary relation $m$ times, yielding an $(n+m)$-ary relation.

Define $I_{2,n}(R)$ and $I^1_{2,n}(R)$ as $I_2(R) \cap V^{n+1}$ [one effect of the intersection with $V^{n+1}$ is to restrict second projections of tuples in $I_2(R)$  to appropriately iterated singletons], and define $I^{m+1}_{2,n}(R)$ as $I_{2
,n+m}(I^m_{2,n}(R))$:   this takes
into account the fact that the tuples get longer.  The intention here is that $I^m_{2,n}$ represents the result of applying the insertion two operation to an $n$-ary relation $m$ times, yielding an $(n+m)$-ary relation.

\item[type lowering axiom:]  $(\forall x:(\exists y:(\forall z:z \in y \leftrightarrow (\forall w:(w,\{z\}) \in x))))$

\item[definition:]
We define ${\tt TL}(x)$ by $(\forall z:z \in {\tt TL}(x) \leftrightarrow (\forall w:(w,\{z\}) \in x))$.  This is a very strange operation!

We define $\iota^{-1}``x$ as ${\tt TL}(V \times x)$.  This is the set of all elements of singletons
belonging to $x$.  We can then define $\iota``x$, the elementwise image of $x$ under the singleton operation, as $\iota^{-1}``(\iota^2``(x))$.

Further, we define $\iota^{-(n+1)}``(x)$ as $\iota^{-1}``(\iota^{-n}``(x))$ for each concrete natural number $n$, and $\iota^{n+1}``(x)$ as $\iota``(\iota^{n}``(x))$

We develop an important operation step by step.

${\tt TL}(R) = \{z:(\forall w:(w,\{z\})\in R)\}$

Dually, $({\tt TL}(R^c))^c = \{z:(\exists w:(w,\{z\})\in R)\}$

Now $({\tt TL}((R^\iota)^c))^c = \{z:(\exists w:(w,\{z\})\in R^\iota)\}$,

which is the same as $({\tt TL}((R^\iota)^c))^c = \{z:(\exists w:(\{w\},\{z\})\in R^\iota)\}$ because all elements
of the domain of $R^{\iota}$ are singletons,

which is the same as $({\tt TL}((R^\iota)^c))^c = \{z:(\exists w:(w,z)\in R)\}$

so we define ${\tt rng}(R)$ as $({\tt TL}((R^\iota)^c))^c$,
and define ${\tt dom}(R)$ as ${\tt rng}(R^{-1})$.

\item[subset axiom:]  $(\exists x:(\forall yz:(y,z) \in x \leftrightarrow (\forall w:w \in y \rightarrow w \in z)))$

We define $[\subseteq]$ as the intersection of a set provided by the previous axiom with $V \times V$:  $[\subseteq]$ is the set of all ordered pairs $(x,y)$ such that $x \subseteq y$.

We define $[\in]$ as $[\subseteq]\cap (1 \times V)$.

\end{description}

This is not our favorite finite axiomatization of NF (or our favorite finite template axiomatization of TST) but it is the one verified in the Lean formalization and also basically the oldest one, so we present a verification of it in outline at least.

What we need to do is verify that $\{x:\phi\}$ exists for each formula $\phi$ of the language of TST, to ensure that comprehension holds.  We do this by induction on the structure of formulas.

$\{x:\neg \phi\}$ is $\{x :\phi\}^c$.

$\{x :\phi \wedge \psi\}$ is $\{x:\phi \} \cap \{x:\psi\}$.

Now we have the much more complex task of analyzing $$\{t_{n}:(\forall t_i:\phi(t_1,\ldots,t_n))\}.$$

Choose a type $\tau'$ higher than the type $\tau_i$ of each $t_i$.  Do this for bound variables
as well, and further, in each formula $(\forall t_i:\psi)$ or $(\exists t_i:\psi)$ we require that all occurrences of $t_i$ be free in $\psi$ and the index of $t_i$ be less than the index
of any other variable appearing free in $\psi$.  It should be clear that we can do this without loss of generality.

Where the type of $t_i$ is $\tau_i$, we define $x_i$ as $\iota^{\tau'-\tau_i}(t_i)$:  we construct $$\{t_{n}:(\forall t_i:\phi(t_1,\ldots,t_n))\}$$ by defining manipulations which allow us to build sets $\{(x_1,\ldots,x_n):\phi^*(x_1,\ldots,x_n)\}$ in which
all the variables are of the same type.  We write $\phi^*$ to suggest that the formula $\phi$ must be transformed to effect our change of variables:  $t_i = t_j$ is equivalent to $x_i = x_j$, and $t_i \in t_j$ is equivalent to $(x_i,x_j) \in [\in]^{\iota^{\tau'-\tau_j}}$  (the reader will see that we use this representation below, though embedded in larger tuples).  A quantifier over $t_i$ is equivalent to a quantifier over $x_i$ restricted to $\iota^{\tau'-\tau_i}``V$.

For any such representation, we have a type signature $$\iota^{\tau'-\tau_1}``V \times \ldots \times \iota^{\tau'-\tau_n}``V.$$  We abbreviate this as $\tau^*$.

The set $\{(x_1,\ldots,x_n):x_1 = x_n\}$ is obtained as $I_{2,2}^{n-2}(\iota^2``V)\cap \tau^*$.

The set $\{(x_1,\dots,x_n):x_1 = x_i\}$ ($i<n$) is obtained as $I_{2,3}^{i-2}(I_3(\iota^2``V))\cap \tau^*$.

We then can represent the
set $\{(x_1,\ldots,x_n):x_i = x_j\}$ (wlog $i<j$) as \newline $I_{1,n-i+1}^{i-1}(\{(x_1,\ldots,x_{n-i+1}):x_1 = x_{j-i+1}\})\cap \tau^*$.

The set $\{(\{x_1\},x_2):x_1 \in x_2\}$ has already been defined above as $[\in]$.

The set $\{(x_1,\ldots,x_n):(\exists uv:x_i = \iota^{k+1}(u) \wedge x_j = \iota^k(v) \wedge u\in v)\}$ is
$${\tt rng}^2(I_3(R^{\iota^{k+2n}}[\in]) \cap \{(x_1,\ldots,x_{n+2}):x_1 = x_{i+2}\} \cap \{(x_1,\ldots,x_{n+2}):x_2=x_{j+2}\}).$$   This
handles all representations of membership statements between $x_i$'s  in the framework we are using.

The set $\{(x_1,\ldots,x_n):(\exists x_1:(x_1,\ldots,x_n)\in R)\}$ is representable as $I_{1,n-1}({\tt rng}(R))$.

This only allows us to quantify over the lowest numbered variable.  This is actually sufficient.  We have renamed bound variables so that all bound variables with different binders are distinct
and in any subformula $(\exists y:\psi)$ $y$ will have lower index than any variable  free in $\psi$.  Then the quantified variable might not be $x_1$, but if it is $x_k$,
we can use the representation $I^{k}_{1,n-k}({\tt rng}^k(R))$:  strip off the first $k-1$ variables, which are dummies already quantified over,
quantify over $x_k$ (stripping off one more variable) and put $k$ dummies back, as it were.

This gives sufficient machinery to handle representations of all sentences in the language of TST (or stratified sentences in the language of NF) in the format we are using.  However, the final representation of  $\{t_n:\phi(t_1,\ldots,t_n)\}$ obtained in this way
will be of the form $R=\{(\iota^{\tau'-\tau_1}(t_1),\ldots,\iota^{\tau'-\tau_n}(t_n)):\phi(t_1.\ldots,t_n)\}$
where some of the $t_i$'s are bound variables (all values of which will appear), and some are parameters.  Let $P$ be the set  of $i$ such that
$t_i$ is a parameter.

The final computation of $\{t_n:\phi(t_1,\ldots,t_n)\}$ will be as
$$\iota^{-(\tau'-\tau_{n})}``\left({\tt rng}^{n-1}\left(R \cap \bigcap_{i \in P}\{(x_1,\ldots,x_n):x_i = \iota^{\tau'-\tau_i}(t_i)\}\right)\right)$$

where a final bit of computation must be exhibited:  $\{(x_1,\ldots,x_n):x_i = t\}$ is realized as
$I^{i-1}_{1,2}(\{\iota^{2(n-i-1)}(t)\}\times V)\cap \tau^*.$

We state for the record that we think this is a bad finite axiomatization of NF, or finite template axiomatization of TST:  we think the treatment of $n$-tuples is terribly difficult to work with.  But it does work, and the second author chose to verify the consistency of TTT by verifying that the finite (template) axiomatization is satisfied  in the structure for the language of TTT that we define, so it merits discussion here.  It should be noted that any finite template axiomatization of TST could be used;  there is no advantage to this one for the formalization, and there is an oddity, because the one impredicative axiom (the axiom of type lowering) does a lot of extra work, since it is also essential for constructing domains and ranges of functions.  In the formalization, the proof of type lowering was in effect divided into the proof of the existence of $\iota^3``{\tt TL}(x)$, which is predicative, and the proof of the existence of set unions of sets of singletons.

The second author reports that she chose this particular finite axiomatization for verification purposes despite its awkwardness because it uses very few defined concepts, and it is therefore clearer to a reader of the Lean formalization that the proofs indeed say what is required.

\newpage

\section{The model description}\label{s:model_description}

In this section we give a complete description of what we claim is a model of tangled type theory.  Our metatheory is some fragment of ZFC.

\subsection{The abstract supertype framework}
\label{ss:supertype_framework}

We first discuss some abstract properties of the types we will construct, and explain the system of ``supertypes".
Types are indexed by a well-ordering $\leq_\tau$ (from which we define a strict well-ordering $<_\tau$ in the obvious way).
We refer to elements of the domain of $\leq_\tau$ as ``type indices".

We first define a system of ``supertypes" (using the same type labels).

For each element $t$ of ${\tt dom}(\leq_\tau)$ we will define a set $\tau^*_t$, called supertype $t$.

If $m$ is the minimal element of the domain of $\leq_\tau$, we choose a set $\tau^*_m$ as supertype $m$.

$\leq_\tau$ and $\tau^*_m$ are the only parameters of the system of supertypes (which is not a model of TTT, but a sort of maximal structure for the language of TTT).

We describe the construction of $\tau^*_t$, assuming that $t \in {\tt dom}(\leq_\tau)$ and $t \neq m$ and for all $u <_\tau t$, we have
defined $\tau^*_u$.

We define $\tau^*_t$ as $$\left\{X \cup \{\{\tau^*_u:u <_\tau t\}\}:X \subseteq \bigcup_{u <_\tau t}\tau^*_u\right\}.$$

An element of supertype $t>_\tau m$ is a subset of the union of all lower types, with $t^+ = \{\tau^*_u:u <_\tau t\}$ added as an element.

Foundation in the metatheory ensures a clean construction here.  An element $x$ of supertype $t>_\tau m$ is always nonempty with $t^+$ as an element.  The set $t^+$ has supertype $u$ as an element for each $u <_\tau t$, so $t^+$ and so $x$  cannot belong to any supertype $u$ with $u <_\tau t$, by foundation.  The labelling element $t^+$ cannot belong to supertype $t$ by foundation, because an element of supertype $t$ must be nonempty and have $t^+$ as an element.  Further, $t^+$ cannot belong to any supertype $v$ with $t <_\tau v$, because any element of $v$ contains $v^+$ as an element which contains supertype $t$ as an element and any element of supertype $t$ contains $t^+$ as an element, so $t^+ \in v$ would violate foundation in the metatheory.

The membership relations of this structure are transparent:  $x \in_{t,u} y$ ($t <_\tau u$) is defined as
$x \in \tau^*_t \wedge y \in \tau^*_u \wedge x \in y$.   Considerations above show that there are no unintended memberships caused by the labelling elements $t^+$, because the labelling elements cannot themselves belong to any supertype.  Note the presence of $\emptyset_t = \{t^+\}$ in supertype $t$, which has no elements of any type $u <_\tau t$ (and is distinct from $\emptyset_v$ for $v \neq t$).

The system of supertypes is certainly not a model of TTT, because it does not satisfy extensionality.  It is easy to construct
many sets in a higher type with the same extension over a given lower type, by modifying the other extensions of the object of higher type.

The system of supertypes does satisfy the comprehension scheme of TTT.  One can use Jensen's method to construct a model of stratified comprehension with no extensionality axiom from the system of supertypes, and stratified comprehension with no extensionality axiom interprets NFU in a manner described by Marcel Crabb\'e in \cite{marcelsf}.

\begin{proposition}[the generality of the system of supertypes]
Any model of TTT (assuming there are any) is isomorphic to a substructure of a system of supertypes.
\end{proposition}
\begin{proof}
Let $M$ be a model of TTT (more generally, any structure for the language of TTT in which each object outside the base type is determined given all of its extensions).  Let $\leq_M$ be the well-ordering on the types of $M$ and let $m$ be the minimal type of $M$.  We will assume as above that $\leq_\tau$ is a well-ordering of type labels $t$ with corresponding actual types $\tau_t$ of $M$:  of course, we could use the actual types of $M$ as type indices, but we preserve generality this way.    We also assume that the sets implementing the types of $M$ are disjoint (it is straightforward to transform a model in which the sets implementing the types are not disjoint to one in which they are, without disturbing its theory, by replacing each $x \in \tau_t$ with $(x,t)$).

We consider the supertype structure generated by ${\leq_\tau}:={\leq_M}$ and $\tau^*_m := \tau_m$.  We indicate how to define an embedding from $M$ into this supertype structure.

Define $I(x) = x$ for $x \in \tau_m = \tau^*_m$.

If we have defined $I$ on each type $u <_\tau t$, we define $I(x)$, for $x \in \tau_t$, as
$$\bigcup_{u <_\tau t} \{I(y):y \in^M_{u,t} x\} \cup \{\{\tau^*_u:u <_\tau t\}\}.$$

It should be clear that as long as $M$ satisfies the condition that an element of any type other than the base type is uniquely determined given all of its extensions in lower types, $I$ is an isomorphism from $M$ to a substructure of the stated system of supertypes.  A model of TTT, in which any one extension of an element of any higher type in a lower type exactly determines the object of higher type, certainly satisfies this condition.  So the problem of constructing a model of TTT is equivalent to the problem of constructing a model of TTT which is a substructure of a supertype system.
\end{proof}

Some advantages of this framework are that the membership relations in TTT are interpreted as subrelations of the membership relation of the metatheory, while the types are sensibly disjoint.

\begin{remark}
Notice that if $\alpha>_\tau \beta$ are supertypes, and $x \in \tau^*_\alpha$, $x \cap \tau^*_\beta$ is the extension of $x$ over supertype $\beta$.  This will be inherited by a scheme of types $\tau_\gamma$ with each $\tau_\gamma \subseteq \tau^*_\gamma$ if an additional condition holds:  for $\alpha>_\tau\beta$, we will have
for $x \in \tau_\alpha$ that $x \cap \tau^*_\beta$ is the extension of $x$ over supertype $\beta$ as already noted:  for it to be the extension over type $\beta$ we need the general condition $x \cap \tau^*_\delta \subseteq \tau_\delta$ for all $\gamma>\delta$ and $x \in \tau_\gamma$.
\end{remark}

\subsection{Preliminaries of the construction}
\label{ss:preliminaries}

Now we introduce the notions of our particular construction in this framework.
\begin{definition}[model parameters]\label{def:model_params}
Let $\lambda$ be a limit ordinal.  Let $\leq_\tau$ be the order on $\lambda \cup \{-1\}$ which has $-1$ as minimal and agrees otherwise with the usual order on $\lambda$.

Let $\kappa$ be an uncountable regular cardinal (that is, an uncountable  regular initial ordinal).  Sets of cardinality $<\kappa$ we call {\em small\/}.  Sets which are not small we may call {\em large\/}.

We make no assumption about the relative sizes of $\kappa$ and $\lambda$.

Let $\mu$ be a strong limit cardinal with $\kappa<\mu$, $\lambda \leq \mu$ and with the cofinality of $\mu$ at least ${\tt max}(\kappa,\lambda)$.
\end{definition}
\begin{remark}\label{rk:example_model_params}
The minimal model parameters are $\lambda = \omega, \kappa = \omega_1, \mu = \beth_{\omega_1}$.
\end{remark}

\begin{definition}[supertypes]\label{def:supertypes}
Let $\tau^*_{-1}=\tau_{-1}$ be $$\{(\nu,\beta,\gamma,\alpha):\nu<\mu \wedge  \beta \in \lambda\cup \{-1\} \wedge \gamma \in \lambda \setminus \{\beta\}\wedge \alpha<\kappa\}.$$  Note that this completes the definition of the supertype structure we are working in as defined in section \ref{ss:supertype_framework}:  we now have a definite reference
for $\tau^*_\alpha$ for $\alpha\in \lambda$.
The important part of this definition is that $\tau_{-1}^*$ has size $\mu$; the precise form of $\tau_{-1}^*$ is chosen to aid with definition \ref{def:f_map}. This type, while important for the model construction, will not be part of our final model of TTT; the types of the final model  will be indexed by $\lambda$.
\end{definition}

Notice that if $\alpha,\beta$ are type indices, $\alpha\in \beta$ is a convenient short way to say $-1 <_\tau \alpha <_\tau \beta$.

\begin{definition}[extended type index]\label{def:extended_type_index}
A nonempty finite subset of $\lambda \cup \{-1\}$ {with minimum element $-1$} may be termed an {\em extended type index}.  If $A$ is a finite subset of $\lambda\cup \{-1\}$  with at least two elements, $A_1$ is defined as $A \setminus \{{\tt min}(A)\}$, and further $A_{n+1}$ is defined as $(A_n)_1$ where this makes sense:  the notation $A_2$ for $(A_1)_1$ sees use.
\end{definition}

\begin{definition}[atoms, litters and near-litters]\label{def:atom_litter_near_litter}
We may refer to elements of $\tau_{-1}$ as {\em atoms\/} from time to time, though they are certainly not atomic in terms of the metatheory.

A {\em litter\/} is a subset of $\tau_{-1}$ of the form $L_{\nu,\beta,\gamma} = \{(\nu,\beta,\gamma,\alpha):\alpha<\kappa\}$.  The litters make up a partition of type $-1$
(which is of size $\mu$) into size $\kappa$ sets.

On each litter $L =  L_{\nu,\beta,\gamma}$ define a well-ordering $\leq_L$:  $(\nu,\beta,\gamma,\alpha) \leq_L (\nu,\beta,\gamma,\alpha')$  iff $\alpha\leq \alpha'$.
The strict well-ordering $<_L$ is defined in the obvious way.  This well-ordering is used in only one place in the paper (definition \ref{def:approx_star}), and its use could easily be avoided, but we find its concreteness appealing.

A {\em near-litter\/} is a subset of $\tau_{-1}$ with small symmetric difference from a litter.  We define $M \sim N$ as $|M \Delta N|<\kappa$, for $M,N$ near-litters:  in English, we say that $M$ is {\em near\/} $N$ iff $M \sim N$.  Note that nearness is an equivalence relation on near-litters.

We define $N^\circ$, for $N$ a near-litter, as the (necessarily unique) litter $L$ such that $L \sim N$.
\end{definition}

In our model, any near-litter will be the extension of $\kappa$-amorphous sets in each type above the base type, and we will use these amorphous sets to limit how much the model knows about the relationship between the different extensions of a model element.

\begin{proposition}\label{def:count_near_litters}
There are exactly $\mu$ near-litters.

\end{proposition}

\begin{proof}
A near-litter is determined as the symmetric difference of a litter $L$ (there are $\mu$ litters) and a small subset of
$\tau_{-1}$.
So it is sufficient to show that a set of size $\mu$ has no more than $\mu$ small subsets.
As $\mu$ has cofinality at least $\kappa$, each small subset of $\mu$ has a strict upper bound $\nu<\mu$, and so is an element of $\bigcup_{\nu < \mu} \mathcal P(\nu)$.
But as $\mu$ is a strong limit cardinal, each $\mathcal P(\nu)$ has size less than $\mu$, so there are only $\mu$ bounded subsets of $\mu$; in particular, there are only $\mu$ small subsets of $\mu$.
\end{proof}

\subsection{Hypotheses for the recursion}\label{ss:hypotheses}

The construction of the types $\tau_\alpha$ is by a recursion.  The initial type $\tau_{-1}$ has already been defined.  For an ordinal $\alpha<\lambda$ we are assuming that $\tau_\beta$ for $-1\leq \beta<\alpha$ have already been constructed and satisfy various hypotheses of the recursion.  We state hypotheses of the recursion in the following subsections in which the construction of $\tau_\alpha$ is described.  We list them here as well, but they are likely best understood where they are encountered in the construction.  Each of these is enforced for $\alpha$ as well:  most of these conditions are explicitly enforced in the course of the construction of $\tau_\alpha$ in this section, but that \ihref{ih:cardinality} holds for $\alpha$ requires an extensive proof in the next section (theorem \ref{thm:count_elements}).

This set of inductive hypotheses is not minimal, but the way it is presented has been chosen to match its actual usage in the following sections.

\begin{enumerate}

\renewcommand{\labelenumi}{(\textbf{I\theenumi})}

\item \label{ih:subset_tau} We assume that for $\gamma<\beta<\alpha$, if $x \in \tau_\beta$, $x \cap \tau^*_\gamma \subseteq \tau_\gamma$.

\item \label{ih:cardinality} We suppose that each $\tau_\beta$ already constructed is of cardinality $\mu$.  Note that we already know that
$\tau_{-1}$ is of cardinality $\mu$.

\item \label{ih:supports} We further intimate that for each $x \in \tau_\gamma$, $-1\leq \gamma<\alpha$, we have defined objects $S$ for which we say that $S$ is a $\gamma$-support of $x$.  [This is introduced before the definition of support is actually given in \ref{def:support}, but there is  no technical obstacle to defining supports earlier if desired].

\item \label{ih:position} We define $\tau_\beta^+$ for any type index $\beta$ as the collection of $(x,S)$ where $x \in \tau_\beta$ and $S$ is a support of $x$.  We assume that we are provided with  a well-ordering $\leq^+_\gamma$ of order type $\mu$ of $\tau_\gamma^+$ ($-1 \leq \gamma <\alpha$), by postulating an injection $\iota^+_\gamma$ from $\tau_\gamma^+$ into $\mu$ (it does not need to be onto) and defining $x \leq^+_\gamma y$ as $\iota^+_\gamma(x) \leq \iota^+_\gamma(y)$.    For any model element $x$ and support $S$ of $x$, we define $\iota^+_*(x,S)$ as $\iota^+_\gamma(x,S)$, where $x \in \tau_\gamma$.    [Hypotheses of the recursion about maps $\iota^+_\gamma$ are stated {in \ihref{ih:position_constraints}}].

\item \label{ih:typed_near_litter} We provide that for every near litter $N$ and every $\beta<\alpha$, there is a unique element $N_\beta$ of $\tau_\beta$ such that $N_\beta \cap \tau_{-1}=N$.

\item \label{ih:extensionality} We further stipulate that extensionality holds for each $\beta\in \alpha$: that is,  for each $\gamma\in \beta$, any $x \in \tau_\beta$ is uniquely determined by $x \cap \tau_\gamma$; $x$ is uniquely determined by $x \cap \tau_{-1}$ only on the additional assumption that $x \cap \tau_{-1}$ is nonempty.

\item \label{ih:pos_typed_near_litter} We assert that $\iota_*(N) \leq \iota_*^+(N_\delta,S)$ will hold for any near litter $N$ and support $S$ of $N_\delta$, $\delta<\alpha$.  This depends on the definition
of the position function $\iota_*$ which appears below at a convenient point but involves no recursion.  \ihref{ih:pos_typed_near_litter} appears in \ihref{ih:position_constraints}, but is mentioned earlier here because it is used earlier in the argument.

\item \label{ih:pre_extensional} We presume that all elements of $\tau_\beta$, $\beta\in \alpha$, are pre-extensional {(definition \ref{def:pre_extensional})}.

\item \label{ih:extensional} We assume that all elements of $\tau_\beta$'s already constructed are extensional {(definition \ref{def:extensional})}.

\item \label{ih:elements_have_supports} We stipulate that all elements of $\tau_\beta$ ($\beta\in\alpha$) have $\beta$-supports.  {Supports are defined in definition \ref{def:support}.}

\item \label{ih:position_constraints} The conditions constraining the choice of functions $\iota^+_\beta$ ($-1 \leq \beta < \alpha$) are

\begin{enumerate}

\item $\iota_*(t) < \iota^+_*(x,S)$ if $(t,A)$ is in the range of $S$ {(supports are functions; see definition \ref{def:support})} and $x$ is not a typed near-litter.

\item $\iota_*(N) \leq \iota^+_*(N_\gamma,S)$ for each near-litter $N$ and support $S$ of $N_\gamma
$ (which appeared earlier as \ihref{ih:pos_typed_near_litter} for reasons stated there; the notation $N_\gamma$ is defined in definition \ref{def:typed_objects}).

\end{enumerate}

We then define $x \leq^+_\beta y$ as $\iota^+_\beta(x) \leq \iota^+_\beta(y)$.  {We will use the position functions $\iota^+_*$ to perform induction along each $\tau_\beta^+$; these constraints on $\iota^+_*$ ensure that objects that a model element depends on (in some suitable sense) are processed before it in the induction.}

\end{enumerate}

\subsection{Machinery for enforcing extensionality in the model}\label{ss:extensionality_machinery}
We describe the mechanism which enforces extensionality in the substructure of this supertype structure that we will build.

The levels of the structure we will define are denoted by $\tau_\alpha$ for \newline $\alpha \in \lambda \cup \{-1\}$.  As we have already noted, $\tau_{-1}=\tau^*_{-1}$ as defined above.

In defining $\tau_\alpha \subseteq \tau^*_\alpha$ for each $\alpha$, we assume that we have already defined $\tau_\beta$ for each $\beta<\alpha$, and that the system of types $\{\tau_\beta:\beta <_\tau \alpha\}$ already defined satisfies various hypotheses which we will discuss as we go [listed at the end of the previous section].
Elements $x$ of $\tau^*_\alpha$ which we consider for membership in $\tau_\alpha$ will have $x \cap \tau^*_\beta \subseteq \tau_\beta$ for $\beta<\alpha$.  We assume that for $\gamma<\beta<\alpha$, if $x \in \tau_\beta$, $x \cap \tau^*_\gamma \subseteq \tau_\gamma$ \ihref{ih:subset_tau}.

We suppose that each $\tau_\beta$ already constructed ($\beta<_\tau \alpha$) is of cardinality $\mu$.  Note that we already know that
$\tau_{-1}$ is of cardinality $\mu$ \ihref{ih:cardinality}.

We further intimate that for each $x \in \tau_\gamma$, $-1\leq \gamma<\alpha$, we have defined objects $S$ for which we say that $S$ is a support of $x$ \ihref{ih:supports}.  The definition of supports will be given in definition \ref{def:support}.  For the moment, we define $\tau_\gamma^+$ as the set of all $(x,S)$ for which $x \in \tau_\gamma$ and $S$ is a support of $x$.  It is a hypothesis of the recursion
that $\tau_\gamma$ is of cardinality $\mu$, from which it follows that  $\tau^+_\gamma$ is of cardinality $\mu$, since there are $\mu$ supports (as we will see when supports are defined).

We also provide a well-ordering $\leq^+_\gamma$ of order type $\mu$ of $\tau_\gamma^+$ ($-1 \leq \gamma <\alpha$), by postulating an injection $\iota^+_\gamma$ from $\tau_\gamma^+$ into $\mu$ (it does not need to be onto) and defining $x \leq^+_\gamma y$ as $\iota^+_\gamma(x) \leq \iota^+_\gamma(y)$ \ihref{ih:position}.   There are some hypotheses of the recursion about maps $\iota^+_\gamma$ which are stated below.   For any model element $x$ and support $S$ of $x$, we define $\iota^+_*(x,S)$ as $\iota^+_\gamma(x,S)$, where $x \in \tau_\gamma$.

We provide that for every near litter $N$ and every $\beta<\alpha$, there is a unique element $N_\beta$ of $\tau_\beta$ such that $N_\beta \cap \tau_{-1}=N$ (\ihref{ih:typed_near_litter}:  we will quite shortly give a precise description of all extensions of this object).

\begin{definition}[typed objects]\label{def:typed_objects}
If $X$ is a subset of $\tau_\gamma$ and $\gamma<\beta$, we define $X_\beta$ as the unique element $Y$ of $\tau_\beta$
such that $Y \cap \tau_\gamma = X$ (if this exists).  Of course, this notation is only usable to the extent that we suppose that extensionality holds.  Notice
that the notation $N_\beta$ is a case of this.

If $N$ is a near-litter we refer to $N_\beta$ as a {\em typed near-litter\/} (of type $\beta$).
\end{definition}

We further stipulate \ihref{ih:extensionality} that extensionality holds for each $\beta\in \alpha$:  that is, for each $\gamma\in \beta$, any $x \in \tau_\beta$ is uniquely determined by $x \cap \tau_\gamma$;  $x$ is uniquely determined by $x \cap \tau_{-1}$ only on the additional assumption that $x \cap \tau_{-1}$ is nonempty.

\begin{definition}[position function for near-litters and singletons of atoms]\label{def:pos_atom_near_litter}
We posit a bijection
$\iota_*$ from the set of all near-litters and singletons of atoms to $\mu$ with the following properties:

\begin{enumerate}

\item $\iota_*(L) < \iota_*(\{a\})$ if $a \in L$ and $L$ is a litter.

\item  $\iota_*(N^\circ)<\iota_*(N)$ if $N$ is a near-litter which is not a litter

\item $\iota_*(\{a\}) < \iota_*(N)$ if $a \in N \Delta N^\circ$

\end{enumerate}

\end{definition}

The function may be constructed directly: first choose an ordering of type $\mu$ on litters, then put singletons of atoms directly after the litter they are inside (and before all later litters), then put near-litters after the corresponding litter and all singletons of atoms in the symmetric difference, then map each near-litter or singleton of an atom to its position in the resulting order.
We do not run out of room before finishing the construction because $\mu$ has cofinality at least $\kappa$.  We prove more formally that there are such functions just below.

\begin{lemma}[position function generator]\label{lem:position_generator}
If $X$ is a set of size at most $\mu$ and $f : X \to \mathcal P(\mu)$ is such that each $f(x)$ is bounded, then there is an injection $\iota : X \to \mu$ such that $\iota(x) \not\in f(x)$.
\end{lemma}

\begin{proof}

Choose a well-ordering $<_X$ of $X$ of type at most $\mu$.
At each stage $x$, assume we have already constructed $\iota(y)$ for $y <_X x$.
We then define
$$
\iota(x) = \min \left(\mu
\setminus f(x)
\setminus \{
\iota(y) : y < x \}\right). $$
Note that as $f(x)$ is bounded in $\mu$, and $\mu$ has initial order type, the set $\mu \setminus f(x)$ is of cardinality $\mu$.
Since $\{ \iota(y) : y < x \}$ has cardinality less than $\mu$, the difference must be nonempty, so the definition given always succeeds.

\end{proof}

We now show how to construct a position function for near-litters and singletons of atoms satisfying the given constraints.

\begin{proof}

First, let $\iota_0$ be any injection from the set of litters into $\mu$.
Next, use lemma \ref{lem:position_generator}, with $X$ taken to be the collection of all singletons of atoms, setting $f(\{a\}) = \{ \nu : \nu \leq \iota_0(L) \}$ where $L$ is the unique litter such that $a \in L$, to create $\iota_1$.

Next, use lemma \ref{lem:position_generator} again, with $X$ taken to be the collection of all near-litters, setting
$$f(N) = \{ \nu : \nu \leq \iota_0(N^\circ)  \} \cup \bigcup_{a \in N \Delta N^{\circ}}\{ \nu : \nu \leq \iota_1(\{a\})\},$$
to obtain $\iota_2$.

Finally, we define a map $\iota$ from near-litters and singletons of atoms to $3 \times \mu$ by
$$ g(x) =
\begin{cases}
  (0, \iota_0(x)) & \text{if $x$ is a litter} \\
  (1, \iota_1(x)) & \text{if $x$ is a singleton} \\
  (2, \iota_2(x)) & \text{otherwise}
\end{cases} $$
We can then define $\iota_*(x)$ to be the position of $g(x)$ in the restriction of the lexicographic order on $3 \times \mu$ to the image of $g$ (which is of order type $\mu$).

The resulting function $\iota_*$ is then clearly a bijection, and has the required properties by construction, noting that near-litters that happen to be litters will automatically satisfy the third condition.

\end{proof}

An additional property involving $\iota_*$ which is enforced by inductive hypotheses explained later about the maps $\iota_\beta^+$ is that $\iota_*(N) \leq \iota_*^+(N_\delta,S)$ will hold for any near litter $N$ and support $S$ of $N_\delta$ \ihref{ih:pos_typed_near_litter}.  It must be noted here because it is shortly used.

\begin{definition}[$f$ maps, crucially important]\footnote{The use of model elements with support rather than simply model elements as domain elements of the $f$ maps is a substantial contribution of the second author to the mathematics of the paper, above simply verifying the work of the first author.  The proof could be carried out without this, but it is much easier to present with this refinement.  There are other ways in which the second author has contributed to the mathematics, but this one is especially worthy of note.}\label{def:f_map}
We define for each type index $\beta$ less than $\alpha$ and each ordinal $\gamma$ distinct from $\beta$ a function $f_{\beta,\gamma}$ (whose definition does not actually depend on $\alpha$:  it will be the same at every stage).  $f_{\beta,\gamma}$ is an injection from $\tau_\beta^+$ into the set of litters, with range included in  $\{L_{\nu,\beta,\gamma}:\nu < \mu\}$ to ensure that distinct $f$ maps have disjoint ranges.

When we define $f_{\beta,\gamma}(x)$, we presume that we have already defined it for $y <^+_\beta x$.
We define $f_{\beta,\gamma}(x)$ as $L$, where $\iota_*(L)$ is minimal such that

\begin{enumerate}
\item $L \in  \{L_{\nu,\beta,\gamma}:\nu < \mu\}$,

\item  $\iota^+_*(x) <\iota_*(L)$ [and so for any $N \sim L$, $\iota^+_*(x) <\iota_*(N)$, and for any $z \in L$, $\iota^+_*(x) < \iota_*(\{z\})$],

\item and for any $y<_\beta^+ x$, $f_{\beta,\gamma}(y) \neq L$.

\end{enumerate}

\end{definition}

Note that the ranges of distinct $f$ maps are disjoint sets.

We are now going to describe how extensionality is enforced in our model.
Each model element will have a particular extension, called its {\em distinguished extension\/}, from which all {\em alternate\/} extensions can be calculated using the $f$ maps.

\begin{definition}[pre-extensional]\label{def:pre_extensional}
We define the notion of {\em pre-extensional\/} element of $\tau^*_\beta$ ($-1 <\beta \leq \alpha$).   An element $x$ of $\tau^*_\beta$ is pre-extensional iff there is a $\gamma<\beta$ such that (1) $x \cap \tau^*_\gamma \subseteq \tau_\gamma$, and (2) $\gamma=-1$ if
$x \cap \tau_{-1}$ is nonempty or if any $x \cap \tau_\delta$ ($\delta \in \beta$) is empty,  and (3) for each $\delta \in \beta \setminus \{\gamma\}$, $$x \cap \tau_\delta= \{N_\delta:(\exists a \in x\cap \tau_\gamma:\exists S:N \sim f_{\gamma,\delta}(a,S))\}.$$  We say for any $x \in \tau^*_\beta$ and $\gamma$ with this property that $x \cap \tau_\gamma$ is a {\em distinguished extension\/} of $x$.
\end{definition}

We presume that all elements of $\tau_\beta$, $\beta\in \alpha$, are pre-extensional \ihref{ih:pre_extensional}.

Note that we now know how to compute all other extensions of typed near-litters $N_\beta$, because the $-1$-extension of $N_\beta$ is distinguished and this indicates how to compute all the other extensions.

\begin{definition}[alternate extension map]\label{def:a_map}
For any $\delta \in \alpha$ and nonempty subset $a$ of type $\gamma \neq \delta$, we define $A_\delta(a)$ as
$$\{N_\delta:(\exists x \in a:\exists S:N \sim f_{\gamma,\delta}(x,S))\}.$$
For any nonempty subset $x$ of type $\delta$ there is at most one subset $y$ of any type such that $A_\delta(y)=x$.  There cannot be more than one such $y$ in any given type because the $f$ maps are injective.  There cannot be more than one such $y$ in different types because the ranges of $f$ maps with distinct index pairs are disjoint.   We use the notation $A^{-1}(x)$ for this set if it exists, defining a very partial function $A^{-1}$ on nonempty subsets of types.
\end{definition}

Note that any distinguished extension of any type element $x$ is the image under $A^{-1}$ of its other extensions.

Using properties of the $f$ maps, we can prove the following useful well-foundedness constraint.

\begin{proposition}\label{prop:a_map_well_founded}
No nonempty subset of a type has infinitely many iterated images under $A^{-1}$, and there are no cycles in $A^{-1}$.
\end{proposition}
\begin{proof}
Let $a$ be a nonempty subset of type $\gamma$ for which $A^{-1}(a)$ exists.

Since $A^{-1}(a)$ exists, every element of $a$ is of the form $N_\gamma$ where $N$ is a near-litter.  Choose $N_\gamma \in a$ such that $\iota_*(N)$ is minimal.  Note that $N$ is in fact a litter by the constraints in definition \ref{def:pos_atom_near_litter}.

Let $b = A^{-1}(a)$:  we have $a = A_\gamma(b)$, where $b$ is a subset of some $\tau_\delta$ with $\delta \neq \gamma$.  In particular, $N = f_{\delta,\gamma}(u,U)$ for
some $u \in \tau_\delta$ and $U$ a support of $u$.  If $u=M_\delta$ we can further state that $\iota_*(N) > \iota^+_*(M_\delta,U) \geq \iota_*(M)$:  if $b$ itself has an image under $A^{-1}$, the minimum value of ordinals $\iota_*(M)$ for $M_\delta \in b$ will be less than the minimum value of ordinals $\iota_*(N)$ for $N_\gamma\in a$, which establishes that there is an ordinal parameter determined by a nonempty subset of a type which decreases strictly when $A^{-1}$ is applied (if it is applicable), and so no nonempty subset of a type
can have infinitely many iterated images under $A^{-1}$, nor can there be any cycles in $A^{-1}$.
\end{proof}

We now verify that the order conditions in the definition of $f_{\beta,\gamma}$ ensure that distinguished extensions are unique.
\begin{proposition}\label{prop:distinguished_extension_unique}
For any $x \in \tau_\beta$ there is only one set $x \cap \tau_\gamma$ which is a distinguished extension of $x$.
\end{proposition}
\begin{proof}
If any $x \cap \tau_\gamma$ ($\gamma \in \beta$) is empty or if $x \cap \tau_{-1}$ is not empty, $x \cap \tau_{-1}$ is the unique distinguished extension
(if it is empty of course it coincides with all the other extensions).

So, what remains is the case of $x$ with $x \cap \tau_{-1}$ empty and each $x \cap \tau_\gamma$ nonempty for $\gamma<\beta$.
Suppose that $x \cap \tau_\gamma$ and $x \cap \tau_\delta$ were both distinguished.
Then we would have
$$ x \cap \tau_\gamma = A^{-1}(x \cap \tau_\delta) = A^{-1}(A^{-1}(x \cap \tau_\gamma)), $$
which contradicts proposition \ref{prop:a_map_well_founded} as this would imply that $A^{-1}$ has a cycle of length 2.
\end{proof}

\begin{definition}[extensional]\label{def:extensional}
We say that an element of a type is {\em extensional\/} iff
it is pre-extensional and its distinguished extension has an even number of iterated images under $A^{-1}$.
This implies that each of its other extensions has an odd number of iterated images under $A^{-1}$.\footnote{We do know that we are carefully, explicitly, spelling out a construction which looks very much
like the construction of the bijection in the Cantor--Schr\"oder--Bernstein theorem.  But the details of the maps involved are used, so everything must be spelled out.}
\end{definition}

\begin{proposition}[extensionality]\label{prop:extensionality}
Two extensional model elements with any common extension (over a type other than $\tau_{-1}$) are equal.
\end{proposition}
\begin{proof}
If two extensional model elements have an empty extension (over a type other than $\tau_{-1}$) in common, they both have all extensions empty and are equal.  If two extensional model elements have a nonempty extension in common, it will be the distinguished extension of both, or a non-distinguished extension of both, since distinguished and non-distinguished extensions are taken from disjoint classes of subsets of types (when nonempty).
In either case we deduce that the two elements have the same distinguished extension and thus have all extensions the same and are equal.  Note that this gives weak extensionality over $\tau_{-1}$ (many objects have empty extension over type $-1$) but it gives full extensionality over any other type.
\end{proof}

We assume that all elements of $\tau_\beta$'s already constructed are extensional \ihref{ih:extensional}.  This completes the mechanism for enforcement of extensionality in the structure we are defining.

\subsection{Allowable permutations and supports}

A crucial aspect of this is that we will need to define $\tau_\alpha$ so that it has cardinality $\mu$ for the process to continue {\ihref{ih:cardinality}}.  It is certainly not a sufficient restriction to require elements of $\tau_\alpha$ to be extensional:  we will require a further symmetry condition.

We define classes of permutations of our structures.
\begin{definition}[structural permutation]\label{def:structural_permutation}
A {\em $-1$-structural permutation\/} is a permutation of $\tau_{-1}^* = \tau_{-1}$.

A {\em $\beta$-structural permutation\/} ($-1 < \beta \leq \alpha$) is a permutation $\pi$ of $\tau_\beta^*$ such that for each type $\gamma<\beta$ there is a $\gamma$-structural permutation
$\pi_\gamma$ such that $\pi(x) \cap \tau^*_\gamma = \pi_\gamma``(x \cap \tau^*_\gamma)$ for any $x \in \tau^*_\beta$.

The maps $\pi_\gamma$ are referred to as {\em derivatives\/} of $\pi$.  If $\pi$ is a $-1$-structural permutation, $\pi_{-1}$ may be taken to denote $\pi$ itself.

More generally, for any finite subset $A$ of $\lambda \cup \{-1\}$ with maximum $\alpha$ and for any $\alpha$-structural permutation $\pi$,
define $\pi_A$ as $(\pi_{A \setminus \{{\tt min}(A)\}})_{{\tt min}(A)}$ and $\pi_{\{\alpha\}} = \pi$.  The maps $\pi_A$ may be referred to as iterated derivatives of $\pi$.\footnote{There is a silly notational point here:  we might want to suppose $\pi_A$ and $\pi_\alpha$ to be essentially distinguished in some way we do not actually implement (for example by type face) in order to prevent confusion of $\pi_{\{0\}}$ with $\pi_1$.  A similar problem exists due to the identification of finite ordinals with finite sets of ordinals.  However, it can also be noted that $\pi_{\{0,\ldots,n\}}$ and $\pi_{n+1}$ cannot make sense for the same allowable permutation $\pi$, so we think this is harmless.}  It should be clear that a structural permutation is exactly determined by its iterated derivatives which are $-1$-structural.

Where $\pi$ is a $\beta$-structural permutation with $\beta>-1$, we define $\pi^+$ as $\pi_{-1}$.  This is occasionally useful to reduce notational clutter.
\end{definition}

Structural permutations are defined on the supertype structure generally.  We need a subclass of structural permutations which respects our extensionality requirements.

\begin{definition}[allowable permutation]\label{def:allowable_permutation}
A {\em $-1$-allowable permutation\/} is a permutation $\pi$ of $\tau_{-1}$ such that for any near-litter $N$, $\pi``N$ is a near-litter.

A {\em $\beta$-allowable permutation\/} ($\beta \leq \alpha$) is a $\beta$-structural permutation, each of whose derivatives $\pi_\gamma$ is a $\gamma$-allowable permutation (and satisfies the condition that $\pi_\gamma``\tau_\gamma = \tau_\gamma$) and which satisfies a coherence condition relating the $f$ maps and derivatives of the permutation:  for suitable $\gamma,\delta<\beta$, $$f_{\gamma,\delta}(\pi_\gamma(x),\pi_{\gamma}[S]) \sim \pi_\delta^+``f_{\gamma,\delta}(x,S).$$  (where the action of allowable permutations on supports will be defined shortly).
This coherence condition is motivated in remark \ref{rk:motivate_coherence_condition}.
\end{definition}

Note that a $\beta$-allowable permutation is actually defined on the entire supertype structure, though what interests us about it is its actions on objects in our purported TTT model.

We will now define a {\em support\/} to be a particular kind of object that can be used to track how an allowable permutation acts along different derivatives on a small part of our model.

\begin{definition}[support condition]\label{def:support_condition}

A {\em $\beta$-support condition\/} ($-1 \leq \beta \leq \alpha$) is defined as a pair $(x,A)$, where
\begin{enumerate}

\item $A$ is an extended type index with maximum $\beta$.  (Recall that an extended type index has minimum $-1$).

\item and $x\subseteq \tau_{-1}$ is either a singleton or a near-litter, and must be a singleton in the case $\beta=-1$.

\end{enumerate}
\end{definition}

\begin{definition}[support]\label{def:support}
Where $-1\leq\beta \leq \alpha$, a {\em $\beta$-support\/} is defined as a function $S$ from a small ordinal to $\beta$-support conditions.

We may write $S_\delta$ instead of $S(\delta)$.
\end{definition}

\begin{remark}
$-1$-supports behave somewhat differently from $\beta$-supports for $-1 < \beta$, and we consider the notion of $-1$-support to be merely a technical convenience, often making the statements of definitions and lemmas more uniform.
The condition that $-1$-support conditions cannot contain near-litters is used exactly once, in the proof of proposition \ref{prop:unions_of_singletons}.
\end{remark}

\begin{definition}[operations on supports]\label{def:support_operations}
We define various operations to manipulate supports.
\begin{enumerate}
\item For any support condition $(x,B)$ we define $(x,B)^{\uparrow A}$ as $(x,B\cup A)$ if all elements of the set $A$ dominate all elements of the set $B$.
Further, if $S$ is a support, we define $S^{\uparrow A}$ so that $(S^{\uparrow A})_\epsilon = (S_\epsilon)^{\uparrow A}$.  By an abuse of notation we may write $(x,B)^{\uparrow \beta}$ or $S^{\uparrow \beta}$ where $\beta$ is an ordinal for $(x,B)^{\uparrow \{\beta\}}$ or $S^{\uparrow \{\beta\}}$.

\item For any supports $S$ and $T$ we denote by $S+T$ a support which consists
of $S$, followed by $T$:  what this means is that $(S+T)_\epsilon = S(\epsilon)$ [which we write $S_\epsilon$] for $\epsilon$ in the domain of $S$, $(S+T)_{{\tt dom}(S)+\epsilon} = T_\epsilon$ for $\epsilon$ in the domain of $T$.

\item We define the action of a $\beta$-allowable permutation $\pi$ on a $\beta$-support $S$:  if $S(\delta) = (x,A)$, $\pi[S](\delta) = (\pi_A``x,A)$.
{This means that $\pi[S] = S$ iff $\pi_A``x = x$ for all $(x,A) \in {\tt rng}(S)$.}

\item
An element $x$ of $\tau^*_\beta$ {\em has $\beta$-support $S$\/} iff for every $\beta$-allowable permutation $\pi$, if $\pi[S] = S$ then $\pi(x)=x$.  We say that  an element $x$ of $\tau^*_\beta$ which has a $\beta$-support is {\em $\beta$-symmetric.}  We allow ourselves to say briefly ``$x$ has support $S$" and ``$x$ is symmetric" if $x \in \tau_\beta$ and $S$ is a $\beta$-support of $x$.
\end{enumerate}
\end{definition}

\begin{remark}[counting supports]\label{rk:counting_supports}

It is straightforward to observe that there are $\mu$ $\beta$-supports for $\beta\leq \alpha$:  there are $\mu$ atoms, $\mu$ near-litters, and
$<\mu$ finite subsets of $\beta +1<\lambda \leq \mu$ (all type indices involved in a $\beta$-support are $\leq \beta$);  thus the set of $\beta$-support conditions (which we will call {\tt SC} for the moment) is of size $\mu$;
note that the set $\kappa \times {\tt SC}$ is of cardinality $\mu$ and each $\beta$-support is a small subset of  $\kappa \times {\tt SC}$, and so, as we have already seen than sets of size $\mu$ have $\mu$ small subsets, it follows that there are no more than $\mu$ $\beta$-supports.
Thus $\tau_\beta^+$ is already known to be of size $\mu$ for $\beta< \alpha$.

\end{remark}

It is important to note that if $S$ is a support of $x\in \tau_\beta$, $\pi[S]$ is a support of $\pi(x)$ for any $\beta$-allowable permutation $\pi$.

\begin{remark}[motivation of the coherence condition in definition \ref{def:allowable_permutation}]\label{rk:motivate_coherence_condition}
Recall the coherence condition for $\beta$-allowable permutations ($\beta \leq\alpha$): for any $(x, S) \in \tau_\gamma^+$ and any $\delta \in \beta \setminus \{\gamma\}$,
$$f_{\gamma,\delta}(\pi_\gamma(x),\pi_{\gamma}[S]) \sim \pi_\delta^+``f_{\gamma,\delta}(x,S).$$
The motivation for this is that we need $\beta$-allowable permutations to send extensional elements of supertypes to extensional elements.  Suppose $x \in \tau_\beta$ and
$x \cap \tau_\gamma = \{b\}$.  If $x$ is extensional, this has to be the distinguished extension of $x$.  For any $\delta \in \beta \setminus \{\gamma\}$,
it follows that $x \cap \tau_\delta$ is the set of all $N_\delta$ such that $N \sim f_{\gamma,\delta}(b,S)$ for some support $S$ of $b$.  This tells us that a $\beta$-allowable permutation $\pi$, such that $\pi(x)$ has $\gamma$-extension $\{\pi_\gamma(b)\}$, must have the  $\delta$-extension of $\pi(x)$ equal to $$\pi_\delta``\{N_\delta:\exists S:N \sim f_{\gamma,\delta}(b,S)\}$$
but must also have its $\delta$-extension equal to $$\{N_\delta:\exists S:N \sim f_{\gamma,\delta}(\pi_\gamma(b),S)\}.$$  This tells us that $$\pi_\delta(f_{\gamma,\delta}(b,S)_\delta) \in \{N_\delta:(\exists T:N \sim f_{\gamma,\delta}(\pi_\gamma(b),T))\}$$ for each support $S$ of $b$.  The coherence condition enforces this neatly, showing that it is motivated by considerations required to get extensionality to work: the action of $\pi_\gamma$ conveniently correlates supports of $b$ with supports of $\pi_\gamma(b)$.

\end{remark}

\begin{proposition}[allowable permutations preserve extensionality]\label{prop:allowable_preserves_extensionality}
Allowable permutations map extensional elements of supertypes to extensional elements.
\end{proposition}
\begin{proof}
Recall that for nonempty  $a\subseteq\tau_\gamma$ we defined $A_\delta(a)$ in definition \ref{def:a_map} as $$\{N_\delta:(\exists x \in a:(\exists S:N \sim f_{\gamma,\delta}(a,S)))\}.$$

If $\pi$ is allowable of suitable index, $\pi_\delta``A_\delta(a)= A_\delta(\pi_\gamma``a)$ follows from the coherence condition.  Verify this:

Suppose we have $N_\delta$ with $x \in a$ such that $N \sim f_{\gamma,\delta}(x,S)$.  Then $$\pi_\delta(N_\delta)  \cap \tau_{-1} = (\pi_\delta)_{-1}``N \sim (\pi_\delta)_{-1}``f_{\gamma,\delta}(x,S) \sim f_{\gamma,\delta}(\pi_\gamma(x),\pi_\gamma[S]).$$  So any element of $\pi_\delta``A_\delta(a)$ is in $A_\delta(\pi_\gamma``a)$.

Suppose we have $N_\delta$ with $x \in a$ such that $N \sim f_{\gamma,\delta}(\pi_\gamma(x),S)$.  We then have $N \sim (\pi_\delta)_{-1}``f_{\gamma,\delta}(x,\pi_\gamma^{-1}[S])$.  We want to show that $\pi_\delta^{-1}(N_\delta) \in A_\delta(a)$.  We have $$\pi_\delta^{-1}(N_\delta) \cap \tau_{-1} = (\pi_\delta)_{-1}^{-1}``N \sim
(\pi_\delta)_{-1}^{-1}``((\pi_\delta)_{-1}``f_{\gamma,\delta}(x,\pi_\gamma^{-1}[S])) = f_{\gamma,\delta}(x,\pi_\gamma^{-1}[S]),$$ establishing what we need.

Notice that this shows that the coherence condition implies that the image under an allowable permutation of a pre-extensional element of our structure is pre-extensional.

Now this implies that if $a \subseteq \tau_\gamma$, then $A^{-1}(a)$ exists and is in $\tau_\delta$ exactly if $A^{-1}(\pi_\gamma``a)$ exists and is in $\tau_\delta$, and moreover $A^{-1}(\pi_\gamma``a)$ is equal to $\pi_\delta``A^{-1}(a)$ if it exists under these conditions.  This verifies that the coherence condition implies that allowable permutations preserve full extensionality, not just pre-extensionality:  the number of iterated images under $A^{-1}$ of an extension that exist is not affected by application of an allowable permutation in a suitable sense.
\end{proof}

\subsection{Model elements defined}\label{ss:model_elements}

\begin{definition}[model elements]
We stipulate that all elements of $\tau_\beta$ ($\beta\in \alpha$) have $\beta$-supports [enforcing \ihref{ih:elements_have_supports}], and define $\tau_\alpha$ as the set of elements $x$ of $\tau^*_\alpha$ such that
$x \cap \tau^*_{\beta} \subseteq \tau_\beta$ for each $\beta<\alpha$, $x$ is extensional, and $x$ has an $\alpha$-support.

\end{definition}

Note that an image of an element of $\tau_\beta$ ($\beta\leq \alpha$)  under a $\beta$-allowable permutation will belong to $\tau_\beta$, because supportedness and extensionality are preserved by allowable permutations.

The definition explicitly enforces \ihref{ih:subset_tau}, \ihref{ih:supports}, \ihref{ih:typed_near_litter} (a near-litter obviously has a support), \ihref{ih:extensionality}, \ihref{ih:pre_extensional}, \ihref{ih:extensional}, \ihref{ih:elements_have_supports} for subsequent stages of the construction.

We still have to prove that the cardinality of $\tau_\alpha$, and so of $\tau^+_\alpha$, is $\mu$, to show that the construction works (verification of \ihref{ih:cardinality} for subsequent stages is in the next section).
{However, we can show now that given \ihref{ih:cardinality}, we can satisfy the remaining hypotheses \ihref{ih:position}, \ihref{ih:pos_typed_near_litter}, \ihref{ih:position_constraints}.}

It should be noted that type 0 has a very simple description:  the $-1$-extensions of type 0 objects are exactly the sets with small symmetric difference from small or co-small unions of litters, and that these are the same extensions over type $-1$ which appear in any positive type.

\begin{proposition}[symmetric comprehension]\label{prop:symmetric_comprehension}
Let $X$ be a subset of a type $-1 < \beta < \alpha$.
Suppose that there is an $\alpha$-support $S$ with the property that for any $\alpha$-allowable $\pi$, when $\pi[S] = S$, we have $\pi_\beta``X = X$.
Then $X = x \cap \tau_\beta$ for some (necessarily unique) model element $x \in \tau_\alpha$.
\end{proposition}
\begin{proof}
If $X$ is empty, then $X$ is simply the $\beta$-extension of the empty model element.
Suppose $X$ is nonempty, and let $n$ be the number of iterated images that $X$ has under $A^{-1}$.

If $n$ is even, consider the element $x \in \tau^*_\alpha$ such that $x \cap \tau^*_\beta = X$, $x \cap \tau^*_{-1} = \emptyset$, and whose other extensions at levels $-1 < \gamma < \alpha$, $\gamma \neq \beta$ are given by $x \cap \tau^*_\gamma = A_\gamma(X)$.
We show this is a model element.
Clearly each extension satisfies $x \cap \tau^*_\gamma \subseteq \tau_\gamma$.
It is pre-extensional because $x \cap \tau_\beta$ is its distinguished extension, and it is extensional as $n$ is even.
It has support $S$: if $\pi$ is $\alpha$-allowable and fixes $S$, then $\pi_\beta `` (x \cap \tau_\beta) = x \cap \tau_\beta$ by assumption, so by extensionality, $\pi$ must fix $x$.
So $x \in \tau_\alpha$ and $x \cap \tau_\beta = X$.

If $n$ is odd, let $Y = A^{-1}(X) \subseteq \tau_\gamma$.
We claim that $S$ has the property that if $\pi$ is $\alpha$-allowable and $\pi[S] = S$, then $\pi_\gamma``Y = Y$.
Indeed, as allowable permutations commute with $A$ maps by proposition \ref{prop:allowable_preserves_extensionality},
$$\pi_\gamma``Y = \pi_\gamma``A^{-1}(X) = A^{-1}(\pi_\beta``X) = A^{-1}(X) = Y.$$
This allows us to repeat the argument of the previous paragraph using $Y$ to obtain a model element $x \in \tau_\alpha$ satisfying $x \cap \tau_\gamma = Y$, where $x \cap \tau_\gamma$ is distinguished.
Then, by pre-extensionality, we must have $x \cap \tau_\beta = A_\beta(x \cap \tau_\gamma) = X$, as required.
\end{proof}
\begin{remark}
  All model elements for $\alpha > 0$ are of this form: if $x \in \tau_\alpha$ and $\beta \in \alpha$, then any support $S$ for $x$ satisfies the hypothesis that $\pi[S] = S$ implies $\pi_\beta``(x \cap \tau_\beta) = x \cap \tau_\beta$.
\end{remark}
\begin{corollary}
The singleton $\{x\}_\alpha$ exists for all $x \in \tau_\beta$.
\end{corollary}
\begin{proof}
If $S$ is a support for $x$, then $S^{\uparrow \alpha}$ satisfies the hypothesis of symmetric comprehension as required.
\end{proof}

\begin{corollary}[$\kappa$-completeness of the structure]
Let $\beta \in \alpha$.
For any small subset $X$ of type $\alpha$, there is a model element $y \in \tau_\alpha$ such that $$y \cap \tau_\beta = \bigcup_{x \in X} x \cap \tau_\beta.$$
\end{corollary}
\begin{proof}
Let $S$ be a support whose range is obtained from the unions of ranges of supports for each $x \in X$.
Then, if $\pi[S] = S$, we have
$$\pi_\beta``\left(\bigcup_{x \in X} x \cap \tau_\beta\right) = \bigcup_{x \in X} \pi_\beta``(x \cap \tau_\beta) = \bigcup_{x \in X} \pi(x) \cap \tau_\beta = \bigcup_{x \in X} x \cap \tau_\beta.$$
The conclusion then follows directly from symmetric comprehension.
\end{proof}

\begin{corollary}
As small unions exist and singletons exist, $X_\alpha$ must also exist for any small set $X \subseteq \tau_\beta$.
\end{corollary}

\begin{proposition}[position functions]
Given that $\tau_\beta^+$ has size $\mu$ \ihref{ih:cardinality}, there exist a choice of injections $\iota^+_\beta$
enforcing \ihref{ih:position}, \ihref{ih:pos_typed_near_litter}, \ihref{ih:position_constraints}.
More precisely,

\begin{enumerate}

\item $\iota_*(t) < \iota^+_\beta(x,S)$ if $(t,A)$ is in the range of $S$ and $x$ is not a typed near-litter \ihref{ih:position_constraints}.

\item $\iota_*(N) \leq \iota^+_\beta(N_\beta,S)$ for any near-litter $N$ and support $S$ of $N$
  \ihref{ih:pos_typed_near_litter}.

\end{enumerate}

We then define $x \leq^+_\beta y$ as $\iota^+_\beta(x) \leq \iota^+_\beta(y)$.
\end{proposition}
\begin{proof}
We use the position function generator (proposition \ref{lem:position_generator}): recall that this states that if $X$ is a set of size at most $\mu$ and $f : X \to \mathcal P(\mu)$ is such that each $f(x)$ is bounded, then there is an injection $\iota : X \to \mu$ such that $\iota(x) \not\in f(x)$.
Setting $X = \tau^+_\beta$, we define
$$ g(x, S) = \bigcup_{N \text{ s.t. } x = N_\delta} \iota_*(N) \cup \bigcup_{(y,A) \in {\tt rng}(S)} \iota_*(y) $$
and let $f(x, S) = \{ \nu : (\exists \xi \in g(x, S) : \nu \leq \xi) \}$.
As $\mu$ has cofinality at least $\kappa$ and each $g(x, S)$ is small, each $f(x, S)$ is bounded as required.
\end{proof}

At this point we have a complete description of the structure which we claim is a model of TTT.

\newpage

\section{Verification that the structure defined is a model}\label{s:verification}

\subsection{The Freedom of Action theorem}\label{ss:foa}

In this section, we will prove a theorem that gives us a way to construct allowable permutations (theorem \ref{thm:foa}).
To do this, we will define a certain class of objects for each type, called {\em approximations\/}, that behave like partial structural permutations.
We will say that an approximation is called {\em coherent\/} if it satisfies a particular version of the coherence condition on allowable permutations (definition \ref{def:allowable_permutation}).
We will then show that every coherent approximation can be identified with the restriction of some allowable permutation.

\begin{definition}[$-1$-approximation]\label{def:base_approx}
  A {\em $-1$-approximation\/} is a function $\psi$ such that:
  \begin{enumerate}
    \item The domain and image of $\psi$ are the same and $\psi$ is injective.
    \item Each domain element of $\psi$ is either an atom or a litter, and moreover, $\psi$ maps atoms to atoms and litters to litters.
    \item For each litter $L$, $\psi$ and $\psi^{-1}$ are defined on only a small collection of atoms $a \in L$.
  \end{enumerate}
\end{definition}
We will associate a partial function $\psi^*$ on atoms to each $-1$-approximation $\psi$.
The action of $\psi$ on atoms will agree with the action of $\psi^*$, and the action of $\psi$ on litters will agree with the pointwise action of $\psi^*$ only up to nearness.

\begin{definition}\label{def:approx_star}
  If $\psi$ is a $-1$-approximation, we define the partial function $\psi^*$ by:
  \begin{enumerate}
    \item If $a$ is an atom and $a \in {\tt dom}(\psi)$, then $\psi^*(a) = \psi(a)$.
    \item If $a$ is an atom with $a \notin {\tt dom}(\psi)$ but $a \in L$ and $L \in {\tt dom}(\psi)$, then
    $$ \psi^*(a) = \pi_{M, N}(a),\ \mathrm{where}\ M = L \setminus {\tt dom}(\psi) \ \mathrm{and}\ N = \psi(L) \setminus {\tt dom}(\psi) $$
    where for any co-small subsets of litters $M, N$, the map $\pi_{M,N}$ is the unique map $\rho$ from $M$ to $N$ such that for any $x, y \in M$,
    $$x <_{M^{\circ}} y \leftrightarrow \rho(x) <_{N^\circ} \rho(y):$$
    $\pi_{M,N}$ is the unique map from $M$ onto $N$ which is strictly increasing in the order determined by fourth projections of atoms.  Notice that $ \pi_{M,N} \circ \pi_{L,M} = \pi_{L,N}$ will hold if $L,M,N$ are co-small subsets of litters,
and $\pi_{L,M} \circ \pi_{M,L} = \pi_{L,L}$ which is the identity map on $L$, under the same conditions.\footnote{The choice of these maps does not need to be so concrete, but the fact that it can be indicates for example that there is no use of choice here.  We like the concreteness of this approach.}
  \end{enumerate}
\end{definition}
\begin{remark}
  If $N$ is a near-litter and $N \subseteq {\tt dom}(\psi^*)$, then $N^\circ \in {\tt dom}(\psi)$, and additionally $\psi^* `` N \sim \psi(N^\circ)$.
  The converse is not true: $N^\circ \in {\tt dom}(\psi)$ does not imply that $N \subseteq {\tt dom}(\psi^*)$, but it does imply that $N^\circ \subseteq {\tt dom}(\psi^*)$.

  Note that if $n$ is any integer, $\psi^n$ is also a $-1$-approximation with the same domain, where we take the convention that $\psi^0$ is the identity map on ${\tt dom}(\psi)$.
  Using the condition $\pi_{M,N} \circ \pi_{L,M} = \pi_{L,N}$, we obtain the equation $(\psi^n)^* = (\psi^*)^n$.
\end{remark}
From now on, we avoid using the action of $\psi$ directly where possible, and instead use $\psi^*$.
\begin{remark}\label{rk:minus_one_approx_allowable}
  $\psi^*$ is a permutation of atoms.
  If it is defined on all of $\tau_{-1}$, it is a $-1$-allowable permutation.
\end{remark}
\begin{definition}[extension]
  We define a partial order on $-1$-approximations by setting $\psi \preceq_{-1} \chi$ if $\psi \subseteq \chi$ and ${\tt dom}(\psi) \cap \tau_{-1} = {\tt dom}(\chi) \cap \tau_{-1}$.
  That is, $\chi$ may define images for more litters than $\psi$, but may not define images for any new atoms.
  If $\psi \preceq_{-1} \chi$, we may call $\chi$ an {\em extension\/} of $\psi$.
  Note that if $\psi \preceq_{-1} \chi$, then $\psi^* \subseteq \chi^*$.
\end{definition}
\begin{definition}[approximation, in general]\label{def:approx}
  If $-1 < \beta$, a $\beta$-approximation is defined as a function with $\{-1\} \cup \beta$ as domain such that $\psi(\gamma)$ is a $\gamma$-approximation for each $\gamma$ in the domain.  We write $\psi_\gamma$ instead of $\psi(\gamma)$.

  We define some operations on approximations.
  \begin{enumerate}
    \item If $A$ is a finite subset of $\lambda\cup \{-1\}$ with maximum element $\beta$, we define the derivative of $\psi$ along $A$ in the following way: $\psi_A = (\psi_{A \setminus \{{\tt min}(A)\}})_{{\tt min}(A)}$ and $\psi_{\{\beta\}} = \psi$.

  If $\beta=-1$, construe $\psi_{-1}$ as $\psi$.

    \item A $\beta$-approximation $\psi$ acts on a support $S$ by
    $$ (\psi^*[S])_\delta = (\psi_A ^*`` x, A) \ \mathrm{where}\ S_\delta = (x, A) $$
    whenever this is defined for each $\delta \in {\tt dom}(S)$.
    \item We define the partial order $\preceq_\beta$ on $\beta$-approximations by defining $\psi \preceq_\beta \chi$ whenever $\psi_\gamma \preceq_\gamma \chi_\gamma$ for all $\gamma < \beta$, and define an {\em extension\/} of a $\beta$-approximation $\psi$ as an approximation $\chi \succeq_\beta \psi$.

  \end{enumerate}
\end{definition}
\begin{definition}[flexibility]\label{def:flexible}
  A near-litter $N$ is {\em  $A$-flexible\/}, where $A$ is an extended type index,  if $|A| \leq 2$ or $N^\circ$ is not in the range of any $f_{\gamma,{\tt min}(A_1)}$
for $-1 \leq \gamma<{\tt min}(A_2)$.
\end{definition}
\begin{definition}[coherent]\label{def:coherent}
  Let $-1 \leq \beta \leq \alpha$.
  We say that a $\beta$-approximation $\psi$ is {\em coherent\/} (c.f.\ the coherence condition on allowable permutations from definition \ref{def:allowable_permutation}) if:
  \begin{enumerate}
    \item If $L$ is $A$-flexible for some $A$, then $\psi_A^*`` L$ is also $A$-flexible whenever this is defined.
    \item If $A$ is an extended type index with maximum $\beta$ and ${\tt min}(A_1) = \gamma < \beta$, and $\delta < {\tt min}(A_2)$ and $(x, S) \in \tau_\delta^+$ are such that $f_{\delta,\gamma}(x, S) \subseteq {\tt dom}(\psi_A^*)$, then there is some $\delta$-allowable permutation $\pi$ such that
    $$ (\psi_{A_2})_\delta^*[S] = \pi[S] $$
    and additionally that
    $$ \psi_A^* `` f_{\delta,\gamma}(x, S) \sim f_{\delta,\gamma}(\pi(x), \pi[S]). $$
    Note that in this case, if $\pi'$ is any other $\delta$-allowable permutation satisfying
    $$ (\psi_{A_2})_\delta^*[S] = \pi'[S] $$
    then we can conclude
    $$ \psi_A^* `` f_{\delta,\gamma}(x, S) \sim f_{\delta,\gamma}(\pi'(x), \pi'[S]) $$
    as $\pi$ and $\pi'$ have actions that coincide on $S$ and therefore at $x$.
  \end{enumerate}
\end{definition}
\begin{remark}\label{rk:power_deriv_coherent}
  If $\psi$ is coherent, then $\psi^n$ is coherent for any integer $n$, and $\psi_\epsilon$ is coherent for any $\epsilon < \beta$.
  Every $-1$-approximation is coherent.
\end{remark}
\begin{definition}[approximate]
  A $-1$-approximation $\psi$ is said to {\em approximate\/} a $-1$-allowable permutation $\pi$ if $\psi^* \subseteq \pi$.
  If $-1 < \beta$, a $\beta$-approximation $\psi$ is said to {\em approximate\/} a $\beta$-allowable permutation $\pi$ if for each $\gamma < \beta$, $\psi_\gamma$ approximates $\pi_\gamma$.
\end{definition}
\begin{remark}\label{rk:preceq_approximates}
  If $\psi \preceq_\beta \chi$ and $\chi$ approximates $\pi$, then $\psi$ approximates $\pi$.
\end{remark}
\begin{definition}[freedom of action]\label{def:foa}
  We say that {\em freedom of action\/} holds at some type index $\beta$ if every coherent $\beta$-approximation $\psi$ approximates some $\beta$-allowable permutation $\pi$.
\end{definition}
\begin{remark}\label{rk:foa_suffices}
  If $\psi$ is a coherent approximation such that $\psi_A$ is defined on all litters (or equivalently, $\psi_A^*$ is defined on all atoms) for all $A$ containing $-1$, then it is easy to see that $\psi$ approximates a unique allowable permutation $\pi$, and $\pi$ is given by $\pi_A= \psi_A^*$.
  So to prove that freedom of action holds at level $\beta$, it suffices by remark \ref{rk:preceq_approximates} to show that every coherent $\beta$-approximation has a coherent extension $\chi$ such that $\chi_A$ is defined on all litters for all $A$ containing $-1$.
\end{remark}
\begin{lemma}\label{lem:foa_flexible}
  Let $\psi$ be a coherent $\beta$-approximation.
  Then $\psi$ admits a coherent extension $\chi$ such that for each extended type index $A$ with maximum $\beta$, $\chi_A$ is defined on all $A$-flexible litters.
\end{lemma}
\begin{proof}
  The construction
  $$ \chi_A = \psi_A \cup \{ (L, L) : L \notin {\tt dom}(\psi_A^*), L\ \mathrm{is}\ A\mathrm{\mbox{-flexible}} \} $$
  suffices.
\end{proof}
\begin{remark}
  This lemma shows in particular that freedom of action holds at type $-1$.
  Indeed, suppose that $\psi$ is a (coherent) $-1$-approximation and $\chi$ is an extension as above.
  By remark \ref{rk:minus_one_approx_allowable}, $\chi^*$ is a $-1$-allowable permutation as all near-litters are $\{-1\}$-flexible, and $\chi$ clearly approximates it.
  But $\psi \preceq_{-1} \chi$, so $\psi$ also approximates $\chi^*$ (remark \ref{rk:preceq_approximates}).
\end{remark}
\begin{lemma}\label{lem:foa_inflexible}
  Let $\psi$ be a coherent $\beta$-approximation.
  Let $A$ be an extended type index with maximum $\beta$, and let $\gamma = {\tt min}(A_1)$ and $\delta < {\tt min}(A_2)$.
  Let $(x, S) \in \tau_\delta^+$ be such that $(\psi_{A_2})_\delta^*[S]$ is defined.

  Then if freedom of action holds at type level $\delta$, there is a coherent extension $\chi \succeq_\beta \psi$ such that $
  \chi_A^* `` f_{\delta,\gamma}(x, S)$ is defined.
\end{lemma}
\begin{proof}
  The $\delta$-approximation $(\psi_{A_2})_\delta$ is coherent (remark \ref{rk:power_deriv_coherent}), so it approximates some $\delta$-allowable permutation $\pi$ by freedom of action.
  In particular, $((\psi_{A_2}^n)_\delta)^*[S] = \pi^n[S]$ for every integer $n$.
  We intend to add the orbit
  $$ f_{\delta,\gamma}(\pi^n(x), \pi^n[S]) \mapsto f_{\delta,\gamma}(\pi^{n+1}(x), \pi^{n+1}[S]) $$
  to $\psi_A^*$, at least up to nearness.

  Suppose that there is some $n$ such that $\psi_A^*$ is defined on $f_{\delta,\gamma}(\pi^n(x), \pi^n[S])$.
  We have $((\psi_{A_2}^n)_\delta)^*[S] = \pi^n[S]$, so $((\psi_{A_2}^{-n})_\delta)^*[\pi^n[S]] = S$.
  Therefore, as $\psi^{-n}$ is coherent (also by remark \ref{rk:power_deriv_coherent}), we obtain
  $$ (\psi^{-n})_A^* `` f_{\delta,\gamma}(\pi^n(x), \pi^n[S]) \sim f_{\delta,\gamma}(x, S). $$
  Therefore,
  $$ f_{\delta,\gamma}(x, S) \subseteq {\tt dom}(\psi_A^*).$$
  So we do not need to extend $\psi$; we are already done.

  Otherwise, we can extend $\psi$ to an approximation $\chi$ given by
  \[ \chi_A = \psi_A \cup \{ (f_{\delta,\gamma}(\pi^n(x), \pi^n[S]), f_{\delta,\gamma}(\pi^{n+1}(x), \pi^{n+1}[S])) : n \in \mathbb Z \} \]
  and \( \chi_B = \psi_B \) for all \( B \neq A \).
  It is easy to see that \( \chi \) is an extension of \( \psi \), and
  $$ \chi_A^* `` f_{\delta,\gamma}(\pi^n(x), \pi^n[S]) \sim f_{\delta,\gamma}(\pi^{n+1}(x), \pi^{n+1}[S]) $$
  for each integer $n$.
  Finally, we must check that this $\chi$ is coherent, but this holds by construction.
\end{proof}
\begin{theorem}[Freedom of Action]\label{thm:foa}
  Freedom of action holds at all type indices $\beta \leq \alpha$.\footnote{It should be noted that the second author has entirely reorganized and rewritten this proof, and indeed the entire section on Freedom of Action.  The underlying mathematics is the same but the presentation is much cleaner.}
\end{theorem}
\begin{proof}
  By induction, we may assume freedom of action holds at all levels $\delta < \beta$.
  Let $\psi$ be a coherent $\beta$-approximation, and use Zorn's lemma to extend $\psi$ to a maximal coherent extension $\chi$; this step uses the fact that coherence is preserved under suprema of chains of approximations.

  Suppose that $\chi_A^*$ is not defined on all litters for some $A$.
  Let $L$ be the litter with minimal position $\iota_*(L)$ such that there is
  an extended type index $A$
  with maximum element $\beta$ such that $L \nsubseteq {\tt dom}(\chi_A^*)$.

  Suppose that $L$ is $A$-flexible.
  By lemma \ref{lem:foa_flexible}, $\chi$ admits a coherent extension $\varphi$ such that $L \subseteq {\tt dom}(\varphi_A^*)$.
  This contradicts maximality of $\chi$.

  Now suppose that $L$ is not $A$-flexible.
  Writing $\gamma = {\tt min}(A_1)$, there are $\delta < {\tt min}(A_2)$ and $(x, S) \in \tau_\delta^+$ such that $L = f_{\delta,\gamma}(x, S)$.

  We claim that $(\chi_{A_2})_\delta^*[S]$ is defined; this will then give a contradiction by lemma \ref{lem:foa_inflexible}.
  To show this, we must prove that if $(t, B) \in {\tt rng}(S)$, then $(\chi_{A_2})_B^*``t$ is defined.
  By definition \ref{def:f_map}, we have $\iota_*^+(x, S) < \iota_*(L)$, and by \ihref{ih:position_constraints}, $\iota_*(t) \leq \iota_*^+(x, S)$.
  If $t$ is a singleton of an atom and $t \subset M$ for a litter $M$, then by definition \ref{def:pos_atom_near_litter}, we obtain $\iota_*(M) < \iota_*(t)$, and therefore that $(\chi_{A_2})_B^*``t$ is defined by minimality of the position of $L$.
  Alternatively, if $t$ is a near-litter, then if $M$ is any litter such that $M \cap t \neq \emptyset$, definition \ref{def:pos_atom_near_litter} again implies that $\iota_*(M) \leq \iota_*(t)$, so $(\chi_{A_2})_B^*``M$ is defined.
  Combining all such litters, we conclude that $(\chi_{A_2})_B^*``t$ is defined.

  Therefore $\chi_A^*$ must be defined on all litters for all $A$.
  By remark \ref{rk:foa_suffices}, this concludes the proof: $\psi$ approximates the allowable permutation $\pi$ given by $\pi_A = \chi_A^*$.
\end{proof}
\begin{remark}
  The use of Zorn's lemma in the previous proof is merely a technical convenience.
  Its use can be excised by instead computing the value of \( \chi_A^* \) at each atom directly, under the inductive hypothesis that its value at each atom at an earlier position was already computed for all \( A \).

  In fact, all of the proofs of this subsection can be phrased in such a way that the axiom of choice is not invoked.
  In particular, the coherent extension defined in \ref{lem:foa_inflexible} is uniquely determined: all choices of \( \pi \) yield the same extension \( \chi \).
  This means that every coherent approximation defined on all flexible litters has a unique coherent extension defined on all litters (but we will not use this fact).
\end{remark}

\subsection{Freedom of Action for incomplete orbits}

In this subsection, we give an alternative form of the Freedom of Action Theorem (\ref{thm:foa}) that is more useful when we do not already have full orbits of atoms or litters for use in approximations.
\begin{definition}[interference]\label{def:interference}
Let $x, y \subseteq \tau_{-1}$.
Their {\em interference\/} is defined to be the union of the small elements of $\{ x \Delta y, x \cap y \}$, which is a small subset of $\tau_{-1}$.
\end{definition}
Thus the interference between two near-litters $M$ and $N$ is either $M \Delta N$ or $M \cap N$, whichever is small.

\begin{definition}[actions]
  A {\em $-1$-action\/} is a function $\xi$ such that:
  \begin{enumerate}
    \item The domain and image of $\xi$ are small, and $\xi$ is injective.
    \item Each domain element of $\xi$ is either an atom or a near-litter, and moreover, $\xi$ maps atoms to atoms and near-litters to near-litters.
    \item If $\xi$ is defined at near-litters $N_1$ and $N_2$, then $N_1^\circ = N_2^\circ$ if and only if $\xi(N_1)^\circ = \xi(N_2)^\circ$.
    This ensures that $\xi$ gives rise to an injective partial function on litters.
    \item Whenever $a, N \in {\tt dom}(\xi)$, we have $a \in N$ if and only if $\xi(a) \in \xi(N)$.
    \item The interference of domain elements of $\xi$ is included in the domain of $\xi$, and the interference of range elements of $\xi$ is included in the range of $\xi$.
  \end{enumerate}
  Note that this definition ensures that the inverse of a $-1$-action is a $-1$-action.
\end{definition}
\begin{definition}
  For $-1 < \beta$, a {\em $\beta$-action\/} is a function with $\{-1\}\cup \beta$ as domain such that $\xi(\gamma)$ is a $\gamma$-action for each $\gamma < \beta$.
  For convenience, if $a$ is an atom we define $\xi_A(\{a\})$ as $\{\xi_A(a)\}$.
\end{definition}
\begin{definition}
  We define derivatives for actions in the same way as we do for approximations.
  Actions can be applied to supports as follows:
  $$ (\xi[S])_\delta = (\xi_A(x), A) \ \mathrm{where}\ S_\delta = (x, A) $$
  As with approximations, this operation is partial.
  An action is called coherent if:
  \begin{enumerate}
    \item If $N$ is $A$-flexible for some $A$, then $\xi_A(N)$ is also $A$-flexible whenever this is defined.
    \item If $A$ is an extended type index with maximum $\beta$ and ${\tt min}(A_1) = \gamma < \beta$, and $\delta < {\tt min}(A_2)$ and $(x, S) \in \tau_\delta^+$ are such that $f_{\delta,\gamma}(x, S) \sim N \in {\tt dom}(\xi_A)$, then there is some $\delta$-allowable permutation $\pi$ such that
    $$ (\xi_{A_2})_\delta[S] = \pi[S] \ \mathrm{and}\ \xi_A(N) \sim f_{\delta,\gamma}(\pi(x), \pi[S]). $$
  \end{enumerate}
  A $\beta$-action $\xi$ is said to approximate a $\beta$-allowable permutation $\pi$ if $\xi_A``x = \pi_A``x$ whenever $x$ is a singleton of an atom or a near-litter, and the left-hand side is defined.
\end{definition}
\begin{theorem}[Freedom of Action for incomplete orbits]\label{thm:foa_actions}
  Every coherent $\beta$-action approximates some $\beta$-allowable permutation $\pi$.
\end{theorem}
\begin{proof}
  Let $\xi$ be a coherent $\beta$-action.
  First, note that we may assume without loss of generality that whenever $N \in {\tt dom}(\xi_A)$, we have
  \begin{equation}
    N \Delta N^\circ \subseteq {\tt dom}(\xi_A) \ \mathrm{and}\ \xi_A(N) \Delta \xi_A(N)^\circ \subseteq {\tt rng}(\xi_A).
    \tag{$*$}
  \end{equation}
  We can extend any $\xi$ to an action $\xi'$ with this property.  If there is an atom in $N^\circ \setminus N \setminus {\tt dom}(\xi_A)$ then neither $N^{\circ}$ nor indeed any $N' \sim N$ containing the atom is in the domain of $\xi$, and we can map each such atom into an unused litter.
  If there is an atom in $N \setminus N^\circ \setminus {\tt dom}(\xi_A)$ then again $N^{\circ}$ is not in the domain of $\xi$ and all $N' \sim N$ which are in the domain of $\xi$ contain the atom, and we can map each such atom to an unused atom which is in $\xi_A(N)^\circ$ and also inside all $\xi_A$-images of near-litters $N' \sim N$.  This handles near-litters $N$ in the domain of $\xi$;  near-litters in the range of $\xi$ are handled dually.
  This modification preserves coherence, as it only adds images of atoms.

  Now, we will turn $\xi$ into a $\beta$-approximation $\psi$.
  To do this, we need to fill in orbits of litters and atoms.

  We first show how to define orbits of flexible litters; the orbits for inflexible litters will appear automatically by coherence.
  For each $A$, let $g_A$ be the partial function given by $g_A(N^\circ) = \xi_A(N)^\circ$ whenever $N$ is $A$-flexible and $\xi_A(N)$ is defined; the third condition on the definition of actions ensures that this map is well-defined and injective.
  We define the restriction of $\psi_A$ to litters to be a small bijective extension of $g_A$, where any new litters in ${\tt dom}(\psi_A)$ are $A$-flexible and disjoint from domain and range elements of $\xi_A$.

  Next, we show how to define orbits of atoms.
  For each $A$, we define the restriction of $\psi_A$ to atoms to be a small bijective extension of the restriction of $\xi_A$ to atoms, with the property that $a \in N$ iff $\psi_A(a) \in \xi_A(N)$ whenever the right-hand side is defined.
  We now explain how to do this.
  Our construction proceeds in $\omega$ stages, where at each stage, we define $\psi_A$ and $\psi_A^{-1}$ at any atoms contained in the domain or range of the function defined so far.
  We describe the construction of images under $\psi_A$; the inverse is dual.

  Suppose that we need to define $\psi_A(a)$, so in particular $a \notin {\tt dom}(\xi_A)$, and that $\xi_A(N)$ is defined for some $N \ni a$.
  By $(*)$, we know that $a \in N^\circ$, and additionally $a \in N'$ for every $N' \sim N$ with $N' \in {\tt dom}(\xi_A)$ because of the interference conditions.
  So we define $\psi_A(a)$ to be any unused atom both in $\xi_A(N)^\circ$ and also inside all $\xi_A(N')$ for $N' \sim N$.
  Such an atom will always exist because there are only a small amount of atoms excluded by each near-litter $N'$, and there are only a small number of such near-litters.

  If $\xi_A(N)$ is not defined for any $N \ni a$, then we are free to choose any unused atom for $\psi_A(a)$.

  We perform the dual construction to define $\psi_A^{-1}$ on required atoms.
  When the construction concludes, the function $\psi_A$ will still be small because $\kappa$ is uncountable and regular.

  The $\beta$-approximation $\psi$ defined as above is obviously coherent because all of the litters $L$ in the domain of $\psi_A$ are $A$-flexible, so by the Freedom of Action Theorem (\ref{thm:foa}), $\psi$ approximates some $\beta$-allowable permutation $\pi$.

  We claim that $\xi$ approximates this $\pi$.
  The case for atoms is obvious, so it suffices to show that for each near-litter $N$ such that $\xi_A(N)$ is defined, we have $\xi_A(N) = \pi_A `` N$.
  By induction we may assume that $\xi_B(M) = \pi_B `` M$ holds for all extended type indices $B$ and near-litters $M$ with $\iota_*(M) < \iota_*(N)$.

  First, we show that $\xi_A(N)^\circ = (\pi_A `` N)^\circ$.
  If $N$ is $A$-flexible, we have $\xi_A(N)^\circ = \psi_A(N^\circ) = (\pi_A `` N)^\circ$ as required.
  Suppose instead that $\gamma = {\tt min}(A_1)$ and $\delta < {\tt min}(A_2)$ are such that $N \sim f_{\delta,\gamma}(x, S)$ for $(x, S) \in \tau_\delta^+$.
  For any $(t, B)$ in the range of $S$, by coherence of $\xi$ we know that $(\xi_{A_2})_B(t)$ is defined, and moreover by inductive hypothesis, that $(\xi_{A_2})_B(t) = (\pi_{A_2})_B``t$.
  But this shows that $\xi_{A_2}$ agrees with $\pi_{A_2}$ on $S$, and so again by coherence, we can conclude that
  $$\xi_A(N) \sim f_{\delta,\gamma}((\pi_{A_2})_\delta(x), (\pi_{A_2})_\delta[S]) \sim \pi_A `` f_{\delta,\gamma}(x, S) \sim \pi_A `` N.$$

  We now conclude the proof by showing that $\xi_A(N) = \pi_A `` N$.
  This proof depends on $(*)$, as well as the fact that the allowable permutation $\pi$ that we constructed in theorem \ref{thm:foa} has the property that for any atom $a \notin {\tt dom}(\psi)$ and any litter $L$, we have $a \in L$ if and only if $\pi(a) \in (\pi``L)^\circ$.

  Suppose that $a \in \xi_A(N)$.
  If $a \in {\tt rng}(\psi_A)$, then $\psi_A^{-1}(a) \in N$ by construction of $\psi_A$, but $\psi_A^{-1}(a) = \pi_A^{-1}(a)$, so $\pi_A^{-1}(a) \in N$.
  Otherwise, $a \in \xi_A(N)^\circ$ by $(*)$, so $\pi_A^{-1}(a) \in N^\circ$, so $\pi_A^{-1}(a) \in N$ by $(*)$.

  Conversely, suppose that $a \in \pi_A `` N$.
  If $a \in {\tt rng}(\psi_A)$, then $\pi_A^{-1}(a) = \psi_A^{-1}(a)$, so $a \in \xi_A(N)$ by construction of $\psi_A$.
  Otherwise, $\pi_A^{-1}(a) \in N^\circ$ by $(*)$, so $a \in (\pi_A `` N)^\circ = \xi_A(N)^\circ$, so by $(*)$ we conclude that $a \in \xi_A(N)$.
\end{proof}

\newpage
\subsection{Counting orbits of supports:  coding functions and specifications}

Now we argue that (given that everything worked out correctly already at lower types) each type $\alpha$ is of size $\mu$, which ensures
that the construction actually succeeds at every type (verification of \ihref{ih:cardinality} for subsequent stages of the construction is thus completed).

\begin{definition}[strong support]\label{def:strong_support}
A $\beta$-support $S$ is called {\em strong\/}
if it satisfies the additional properties that

\begin{enumerate}

\item if $(x,A)$ and $(y,A)$ are in the range of $S$, then $(\{z\},A)$ is in the range of $S$ for all atoms $z$ in the interference of $x$ and $y$ (defined in definition \ref{def:interference}),

\item and if $(x,A)$ is in the range of $S$, $\delta={\tt min}(A_1)$, $\gamma<{\tt min}(A_2)$, and $x^\circ = f_{\gamma,\delta}(y,T)$, then the range of $T^{\uparrow A_2}$ is a subset of the range of $S$:  supports appearing in inverse images under $f$ of litters which are near the first projections of an element of the support have a type-raised copy (mod reindexing) appearing in the support.
\end{enumerate}
\end{definition}

\begin{remark}
It should be evident that if $\pi$ is a $\beta$-allowable permutation and $S$ is a strong $\beta$-support,
$\pi[S]$ is also a strong $\beta$-support.
\end{remark}

\begin{proposition}\label{prop:extend_to_strong_support}
Every $\beta$-support $S$ is an initial segment of some strong support.
\end{proposition}
\begin{proof}
First, for each $(x,A)$ in the range of $S$ where $\delta={\tt min}(A_1)$, $\gamma<{\tt min}(A_2)$ and $x^\circ = f_{\gamma,\delta}(y,T)$, add the range of $T^{\uparrow A_2}$ to the end of $S$.
Repeating this process $\omega$ times, we obtain a new support $S' \supseteq S$ that satisfies the second property of being a strong support; the domain of $S'$ is still a small ordinal as $\kappa$ is uncountable and regular.
Finally, we append to $S$ all pairs $(\{z\},A)$ where $z$ is in the interference of near-litters $M,N$ where $(M,A)$, $(N,A)$ are in the range of $S$.
This produces a support $S'' \supseteq S'$ which is clearly strong.
\end{proof}

\begin{definition}[coding function]\label{def:coding_function}
For any $(x, S) \in \tau_\beta^+$, we can define a function $\chi_{x,S}$ which sends $T=\pi[S]$ to $\pi(x)$ for every $T$ in the orbit of $S$ under
the action of allowable permutations.  We call such functions {\em coding functions\/}.  Note that if $\pi[S]=\pi'[S]$ then $(\pi^{-1}\circ \pi')[S]= S$, so
$(\pi^{-1}\circ \pi')(x)= x$, so $\pi(x)=\pi'(x)$, ensuring that the map $\chi_{x,S}$ for which we gave an implicit definition is well defined.
\end{definition}

We will now define a certain combinatorial object that will characterize the orbit of a strong support under the action of allowable permutations (although we can define it even when the support in question is not strong).

\begin{definition}[specification]
The {\em specification\/} of a $\beta$-support $S$ is the function $S^*$ with the same domain as $S$ given by the following criteria.
We use the notation $S^*_\epsilon$ for $S^*(\epsilon)$.

\begin{enumerate}

\item If $S_\epsilon$ is $(\{x\},A)$, then $S^*_\epsilon$ is given by
$$ S^*_\epsilon = (0, A, \{ \zeta : x \in \pi_1(S_\zeta) \wedge A = \pi_2(S_\zeta) \}). $$
\item If $S_\epsilon$ is $(N,A)$ where $N$ is $A$-flexible, then $S^*_\epsilon$ is given by
$$ S^*_\epsilon = (1, A, \{ \zeta : (\exists N' : S_\zeta = (N',A) \wedge N \sim N') \} ). $$
\item If $S_\epsilon$ is $(N,A)$ where $N$ is not $A$-flexible, so $|A| > 2$, $\gamma = {\tt min}(A_1)$, and $N \sim f_{\delta,\gamma}(x,T)$ with $-1\leq\delta<{\tt min}(A_2)$, then $S^*_\epsilon$ is given by
$$ S^*_\epsilon = (2, A, \chi_{x,T}, \{ ( \zeta, \eta ) : S_\zeta = (T^{\uparrow A_2})_\eta \} ). $$

\end{enumerate}
\end{definition}

\begin{remark}
It should be evident that for any $\beta$-allowable permutation $\pi$ and $\beta$-support $S$, $(\pi[S])^* = S^*$.  What is less evident and our first target result here is that if $S$ is a strong support then any strong $T$ with $T^* = S^*$ is the image of $S$ under the action of an allowable permutation:  the specifications precisely code the orbits in the strong supports under the allowable permutations.
\end{remark}

\begin{lemma}\label{lem:specification_determines_orbit}
The specification of a strong $\beta$-support exactly determines the orbit in the action of $\beta$-allowable permutations on supports to which it belongs:  if two strong $\beta$-supports have the same specification, they are in the same orbit.
\end{lemma}
\begin{proof}
It is straightforward to see that if $S$ is a $\beta$-support and if $\pi$ is a $\beta$-allowable permutation, then $S^* = (\pi[S])^*$.  The relationships between items in the support recorded in the specification are invariant under application of allowable permutations.

It remains to show that if $S$ and $T$ are strong supports, and $S^*=T^*$ is a specification for both, there is an allowable permutation $\pi$ such that $\pi[S]=T$.

We construct $\pi$ using the Freedom of Action Theorem (\ref{thm:foa_actions}).
Define the $\beta$-action $\xi$ in the following way.

If $S_\epsilon = (M,A)$ for $M$ a near-litter, then $T_\epsilon = (N,A)$ for $N$ a near-litter, and we set $\xi_A(M) = N$.
Note that if $S_\delta = (M,A)$ for some $\delta$, the fact that $S$ and $T$ have the same specification implies that $T_\delta = (N,A)$, ensuring that this is a well-defined function.
Similarly, if $S_\delta = (M',A)$ and $M' \sim M$, then $T_\delta = (N',A)$ for some $N' \sim N$.

If we have $S_\epsilon = (\{x\},A)$, we will have $T_\epsilon = (\{y\},A)$ for some $y$, and so we will similarly set $\xi_A(x) = y$.
Note that if $S_\delta = (\{x\},A)$, then $T_\delta = (\{y\},A)$.

As $S$ and $T$ are strong, this description indeed defines a $\beta$-action $\xi$.
We now check that $\xi$ is coherent.
If we assigned $\xi_A(M) = N$, then $(M,A) = S_\epsilon$ and $(N,A) = T_\epsilon$ for some $\epsilon$.
As $S$ and $T$ have the same specification, $M$ is $A$-flexible if and only if $N$ is $A$-flexible.

Now suppose that $\xi_A(M) = N$ where $M^\circ = f_{\delta,\gamma}(x, U)$ and $\gamma = {\tt min}(A_1)$ and $\delta < {\tt min}(A_2)$.
Then $(M,A) = S_\epsilon$ and $(N,A) = T_\epsilon$ for some $\epsilon$.
As $S$ and $T$ have the same specification, there is some $(y,V) \in \tau_\delta^+$ such that $N^\circ = f_{\delta,\gamma}(y,V)$ and
$$ \chi_{x,U} = \chi_{y,V};\quad \{ ( \zeta, \eta ) : S_\zeta = (U^{\uparrow A_2})_\eta \} = \{ ( \zeta, \eta ) : T_\zeta = (V^{\uparrow A_2})_\eta \}. $$
Since $\chi_{x,U} = \chi_{y,V}$, there is a $\delta$-allowable permutation $\pi$ such that $\pi(x) = y$ and $\pi[U] = V$.
So for any $(t,B)$ in the range of $U$, we have $S_\zeta = (t,B)$ if and only if $T_\zeta = (\pi_B(t),B)$.
Therefore, as $S$ and $T$ are strong, we may conclude that $(\xi_{A_2})_\delta[U] = \pi[U]$, and that
$$ \xi_A(M) = N \sim f_{\delta,\gamma}(y,V) = f_{\delta,\gamma}(\pi(x),\pi[U]) $$
as required for coherence.

Finally, by the Freedom of Action Theorem (\ref{thm:foa_actions}), $\xi$ approximates some $\beta$-allowable permutation $\pi$, and in particular, $\pi$ satisfies $\pi[S] = T$, as required.
\end{proof}

\begin{proposition}\label{prop:count_spec}
Suppose that we already know that there are $<\mu$ $\gamma$-coding functions for each $\gamma<\beta$ (which we will be able to assume by inductive hypothesis).
Then there are $<\mu$ specifications of $\beta$-supports for $\beta\leq \alpha$.
\end{proposition}
\begin{proof}
The elements of the range of a $\beta$-specification are built from finite sets of type indices below $\beta$, sets of (pairs of) elements of $\kappa$ (support domain elements), and $\gamma$-coding functions for $\gamma < \beta$.
As $\beta < \lambda \leq \mu$, there are less than $\mu$ possible finite sets of type indices below $\beta$.
Since $\beta < \lambda \leq {\tt cf}(\mu)$, the sum of cardinals each $<\mu$ indexed by ordinals $<\beta$ will be less than $\mu$, so there are less than $\mu$ possible coding functions that can be used here.
Hence the elements of the range of such a specification are taken from a set of size less than $\mu$, and so the cardinality of the set of specifications is bounded by $\nu^{<\kappa}$ where $\nu < \mu$.
Now, $\nu^{<\kappa}$ is the cardinality of the set of functions from some small ordinal to $\nu$, which is less than or equal to the cardinality of the power set of $\nu \times \kappa$, which in turn is less than $\mu$ because $\mu$ is strong limit.
\end{proof}

\begin{proposition}\label{prop:count_support_orbits}
  Suppose that we already know that there are $<\mu$ $\gamma$-coding functions for each $\gamma<\beta$.
  Then there are less than $\mu$ orbits in supports under $\beta$-allowable permutations.
\end{proposition}
\begin{proof}
For any support $S$ that is not necessarily strong, there is a suitable weaker notion of specification: give the specification for a strong support $T$ of which $S$ is an initial segment, together with the domain of $S$.
These weak specifications are not unique, but they do determine the orbit in which $S$ lies in the allowable permutations.
Since there are less than $\mu$ specifications by proposition \ref{prop:count_spec}, and $\kappa < \mu$, we conclude that there are less than $\mu$ orbits in the $\beta$-supports under the allowable permutations.
\end{proof}

\subsection{Types are of size $\mu$ (because there are $<\mu$ coding functions per type)}

The strategy of our argument for the size of the types is to show that that there are $<\mu$ coding functions\ for each type, which implies that there are no more than $\mu$ (and so exactly $\mu$) elements of each type, since every element of a type is obtainable by applying a coding function (of which there are $<\mu$) to a support (of which there are $\mu$).

We will prove this cardinality bound on coding functions by induction on $\beta \leq \alpha$.
It is important to note that this inductive proof occurs internally to a single instance of the main inductive step; that is, $\alpha$ remains fixed, and we do not need to record the fact that the bound on coding functions holds in our main list of inductive hypotheses \ihref{ih:subset_tau}--\ihref{ih:position_constraints}.

This argument will proceed by showing that all coding functions can be written in a particular form, in terms of coding functions at some lower level $-1 \neq \gamma < \beta$.
In particular, this result will depend on the existence of some ordinal $\gamma \in \beta$, so we prove the result for types 0 and $-1$ explicitly here first.

\begin{lemma}\label{lem:count_coding_function_zero}
There are less than $\mu$ coding functions for type 0 and for type $-1$.
\end{lemma}
\begin{proof}
We describe all coding functions for type 0.  The orbit of a strong 0-support in the allowable permutations is determined by the positions in the support  occupied by near-litters, and for each position in the support occupied by a singleton, the positions, if any, of the near-litters in the support  which include it, and for each position in the support occupied by a near-litter, the positions, if any, occupied by near-litters which are near it (this is exactly the information contained in the 0-specification of the support).
There are no more than $2^\kappa$ ways to specify an orbit.  Now for each such support, there is a natural partition of type $-1$ into near-litters, singletons, and a large complement set.
The partition has $\nu<\kappa$ elements, and there are $2^\nu\leq 2^\kappa$ objects of type 0 with this support, each determined by specifying for each compartment in the partition whether it is to be included or excluded from the set.
So there are $2^\nu\leq 2^\kappa$ coding functions for that orbit in the supports, and so for any support (that need not be strong), there are at most $2^\kappa$ coding functions for its orbit, by first extending it to a strong support.   So there are no more than $2^\kappa<\mu$ coding functions over type 0.

Any type $-1$ coding function is associated with a uniquely determined type 0 coding function of a singleton, so there are no more of the former than of the latter.

\end{proof}

With those two cases handled, we return to focus on the case where there is some ordinal $\gamma \in \beta$.

\begin{definition}[combination]\label{def:combination}
Let $\gamma \in \beta$ be ordinals at most $\alpha$.
Let $s$ be a set of $\beta$-coding functions.
We define the {\em combination\/} of $s$ along $\gamma$ to be the function $C_{s,\gamma}$ from $\beta$-supports to $\tau_\beta$ given by
$$ C_{s,\gamma}(S) = \left( \bigcup_{\chi \in s} \bigcup_{T \supseteq S} \chi(T) \cap \tau_\gamma \right)_\beta. $$
That is, for a support $S$, the $\gamma$-extension of the model element $C_{s,\gamma}(S)$ is precisely the union of the $\gamma$-extensions of model elements obtained by applying coding functions from $s$ to supports extending $S$.

That $C_{s,\gamma}(S)$ is always a model element requires proof: by symmetric comprehension (proposition \ref{prop:symmetric_comprehension}), it suffices to show that when a $\beta$-allowable permutation $\pi$ fixes $S$, $\pi_\gamma$ fixes the set
$$ \bigcup_{\chi \in s} \bigcup_{T \supseteq S} \chi(T) \cap \tau_\gamma, $$
but this follows from the fact that for any $\beta$-allowable $\pi$,
\begin{align*}
  \pi_\gamma``\left(\bigcup_{\chi \in s} \bigcup_{T \supseteq S} \chi(T) \cap \tau_\gamma\right)
  &= \bigcup_{\chi \in s} \bigcup_{T \supseteq S} \pi_\gamma``(\chi(T) \cap \tau_\gamma) \\
  &= \bigcup_{\chi \in s} \bigcup_{T \supseteq S} \pi(\chi(T)) \cap \tau_\gamma \\
  &= \bigcup_{\chi \in s} \bigcup_{T \supseteq S} \chi(\pi[T]) \cap \tau_\gamma \\
  &= \bigcup_{\chi \in s} \bigcup_{T \supseteq \pi[S]} \chi(T) \cap \tau_\gamma.
\end{align*}
The same calculation shows that for any orbit of $\beta$-supports $o$ under the action of allowable permutations, the restriction $C_{s,\gamma} \restriction o$ is a coding function.
\end{definition}

Our claim is that every $\beta$-coding function can be written in this form, where the elements of the set $s$ are particular coding functions of singletons of objects of type $\gamma$.  We will make the reference of ``particular" here more precise.

\begin{definition}[raised singleton]\label{def:raised_singleton}
  Let $\gamma \in \beta$.
  A $\beta$-coding function $\chi$ is called a {\em raised $\gamma$-singleton\/} if there is a $\beta$-support $S$ and $(x,U) \in \tau_{\gamma}^+$ such that $\chi = \chi_{\{x\}_\beta,S + U^{\uparrow \beta}}$.
\end{definition}
\begin{lemma}\label{lem:eq_combine}
  Let $\gamma \in \beta$ be ordinals at most $\alpha$, and let $\chi$ be a $\beta$-coding function.
  Then there is a set $s$ of raised $\gamma$-singletons such that $\chi = C_{s,\gamma} \restriction {\tt dom}(\chi)$.
\end{lemma}
\begin{proof}
  Let $\chi = \chi_{x,S}$.
  The set we claim works is
  $$ s = \{ \chi_{\{y\}_\beta,S + U^{\uparrow \beta}} : y \in x \cap \tau_\gamma \wedge U \text{ supports } y \}. $$
  Since $C_{s,\gamma} \restriction {\tt dom}(\chi)$ is a coding function, it suffices to check that $C_{s,\gamma}(S) = x$, which we show by $\gamma$-extensionality.

  If $y \in C_{s,\gamma}(S) \cap \tau_\gamma$, then there is some $z \in x \cap \tau_\gamma$, some support $U$ of $z$, and $T \supseteq S$ such that
  $$ y \in \chi_{\{z\}_\beta,S + U^{\uparrow \beta}}(T). $$
  Hence there is a $\beta$-allowable $\pi$ such that $\pi[S + U^{\uparrow \beta}] = T$; in particular, $\pi[S] = S$ as $T$ is an extension of $S$, so $\pi(x) = x$.
  Therefore,
  $$ \pi_\gamma^{-1}(y) \in \chi_{\{z\}_\beta,S + U^{\uparrow \beta}}(\pi^{-1}[T]) = \chi_{\{z\}_\beta,S + U^{\uparrow \beta}}(S + U^{\uparrow \beta}) = \{z\}_\beta, $$
  giving $\pi_\gamma^{-1}(y) \in x$ and so $y \in x$ as required.

  Conversely, if $y \in x \cap \tau_\gamma$, then for any support $U$ of $y$, we immediately obtain $y \in C_{s,\gamma}(S)$ from the fact that $\chi_{\{y\}_\beta,S + U^{\uparrow \beta}} \in s$.
\end{proof}
\begin{lemma}\label{lem:count_coding_function}
  There are less than $\mu$ coding functions for every type index $\beta \leq \alpha$.
\end{lemma}
\begin{proof}
  We may assume by induction that there are less than $\mu$ coding functions of every type $\gamma < \beta$.
  Since the cases for $\beta = 0, -1$ were proven in lemma \ref{lem:count_coding_function_zero}, we may additionally assume that there exists some ordinal $\gamma \in \beta$.
  By lemma \ref{lem:eq_combine}, every $\beta$-coding function can be written in the form $C_{s,\gamma} \restriction d$, where $s$ is a set of raised $\gamma$-singletons and $d$ is an orbit of $\beta$-supports under the action of allowable permutations.
  By proposition \ref{prop:count_support_orbits}, there are less than $\mu$ such orbits, so it suffices to show there are less than $\mu$ possible sets of raised $\gamma$-singletons; indeed, as $\mu$ is a strong limit, it further suffices to simply show that there are less than $\mu$ raised $\gamma$-singletons.

  We show that any raised singleton $\chi_{\{x\}_\beta,S + U^{\uparrow \beta}}$ is determined by ${\tt dom}(S)$, the orbit of $S + U^{\uparrow \beta}$ under the allowable permutations, and the $\gamma$-coding function $\chi_{x,U}$.
  Indeed, suppose that $S$ and $T$ have the same domain, $S + U^{\uparrow \beta}$ and $T + V^{\uparrow \beta}$ share the same orbit under the $\beta$-allowable permutations, and the coding functions $\chi_{x,U}$ and $\chi_{y,V}$ coincide.
  Then there is a $\beta$-allowable $\pi$ such that $\pi[S] = T$ and $\pi_\gamma[U] = V$, and additionally there is a $\gamma$-allowable $\pi'$ such that $\pi'[U] = V$ and $\pi'(x) = y$.
  Therefore, we additionally have that $\pi_\gamma(x) = y$ as $U$ supports $x$ and both $\pi_\gamma$ and $\pi'$ have the same action on $U$.
  Hence, the two raised singletons must coincide.

  It now simply remains to note that there are less than $\mu$ such orbits of $\beta$-supports under the allowable permutations by proposition \ref{prop:count_support_orbits}, and that we assume by inductive hypothesis that there are less than $\mu$ coding functions of type $\gamma$, in order to conclude the desired result.
\end{proof}

Thus, we conclude
\begin{theorem}\label{thm:count_elements}
Each type $\tau_\alpha$ has exactly $\mu$ elements.
\end{theorem}
\begin{proof}
Any element of a type is determined by a support (of which there are $\mu$ by remark \ref{rk:counting_supports}) and a coding function (there are $<\mu$ of these by lemma \ref{lem:count_coding_function}), so a type has no more than $\mu$ elements (and obviously has at least $\mu$ elements).
\end{proof}

\newpage
\subsection{The structure is a model of predicative TTT}\label{ss:predicative_ttt}

{There is then a very direct proof of the following:}
\begin{proposition}\label{prop:predicative_ttt}
  The structure presented is a model of predicative TTT (in which the definition of a set at a particular type may not mention any higher type).
\end{proposition}
\begin{proof}

Use $E$ for the membership relation $\in_{TTT}$ of the structure defined above (in which the membership of type $\beta$ objects in type $\alpha$ objects is actually a subrelation of the membership relation of the metatheory, a fact inherited from the scheme of supertypes).  It should be evident that $x E y \leftrightarrow \pi_\beta(x) E \pi(y)$,
where $x$ is of type $\beta$, $y$ is of type $\alpha$, and $\pi$ is an $\alpha$-allowable permutation.

Suppose that we are considering the existence of $\{x : \phi^s\}$, where $\phi$ is a formula of the language of TST with $\in$ translated as $E$, and $s$ is a strictly increasing sequence of types\footnote{We use notation $s(i)$ instead of $s_i$ to avoid undue nesting of subscripts, in what follows.}.  Note that the TTT type of $x$ is some $s(j)$ and the TTT type of $\{x : \phi^s\}$ is $s(j+1)$.  For any $s(j+1)$-allowable permutation $\pi$, the truth value of each subformula of $\phi$ will be preserved if we replace each $u$ of type $s(i)$ with $\pi_{A_{s,i}}(u)$, where  $A_{s,i}$ is the set of all $s(k)$ for $i \leq k \leq j+1$ [there being no variables of type higher than $s(j+1)$]:  $\pi_{A_{s,i}}(x) E  \pi_{A_{s,i+1}}(y)$ is equivalent to $(\pi_{A_{s,i+1}})_{s(i)}(x) E \pi_{A_{s,i+1}}(y)$, which is equivalent to $xEy$ by the observation above. The formula $\phi$ will contain various parameters $a_i$ of types $s(n_i)$ and it is then evident that the set $\{x : \phi^s\}$ will be fixed by any $s(j+1)$-allowable permutation $\pi$ such that $\pi_{A_{s,n_i}}$ fixes $a_i$ for each $i$.  But this means that
 $\{x : \phi^s\}_{s(j+1)}$ is symmetric and belongs to type $s(j+1)$:
 we can merge the supports of the $a_i$'s (with suitable raising of indices) into a single $s(j+1)$-support.  Notice that we assumed the predicativity condition that no variable more than one type higher than $x$ appears (in the sense of TST).

This procedure will certainly work if the set definition is predicative (all bound variables are of type no higher than that of $x$, parameters at the type
of the set being defined are allowed), but it also works for some impredicative set definitions.
\end{proof}

There are easier proofs of the consistency of predicative tangled type theory;
there is a reason of course that we have pursued this one.

It should be noted that the construction given here is in a sense a Fraenkel--Mostowski construction, though we have no real need to reference the usual
FM constructions in ZFA here.  Constructions analogous to Fraenkel--Mostowski constructions can be carried out in TST using permutations of type 0;  here we are doing something much more complicated involving many permutations of type $-1$ which intermesh in precisely the right way.  Our explanation of our technique is self-contained, but we do acknowledge this intellectual debt.

\newpage
\subsection{Impredicativity:  verifying the axiom of union}\label{ss:impredicativity}

What remains to complete the proof is that typed versions of the axiom of set union hold.  That this is sufficient is a fact about predicative type theory.
If we have predicative comprehension and union, we note that for any formula $\phi$, $\{\iota^k(x):\phi(x)\}$ will be predicative if $k$ is taken to be large enough, then application of union $k$ times to this set will give $\{x:\phi(x)\}$.  $\iota(x)$ here denotes $\{x\}$.  It is evidently sufficient to prove that unions of sets of singletons exist.
So what we need to show is the following result.

\begin{proposition}\label{prop:unions_of_singletons}
If  $\alpha>\beta>\gamma$ and $G \subseteq \tau_\gamma$, and $$\{\{g\}_\beta:g \in G\}_\alpha$$ is symmetric (has an $\alpha$-support, so belongs to $\tau_\alpha$), then $G_\beta$ is symmetric (has a $\beta$-support, so belongs to $\tau_\beta$).
\end{proposition}
\begin{proof}
Suppose that $\{\{g\}_\beta:g \in G\}_\alpha$ is symmetric.  It then has a strong support $S$.  We claim that $S_{(\beta)}$, defined as $\{(z,C): {\tt max}(C)=\beta \wedge (z,C\cup \{\alpha\}) \in S\}$,  is a $\beta$-support for $G_\beta$, where without loss of generality we assume that $S$ has been reordered such that the domain of $S_{(\beta)}$ (and $(S_{(\beta)})_{(\gamma)}$) is an ordinal.

Any $g \in G$ has a strong $\gamma$-support $T$ which extends $(S_{(\beta)})_{(\gamma)}$.

Suppose that the action of a $\beta$-allowable permutation $\pi$ fixes $S_{(\beta)}$.

Our plan is to use Freedom of Action technology to construct an $\alpha$-allowable permutation $\pi^*$ whose action on $S$ is the identity
and whose action on $T^{\uparrow\{\alpha,\beta\}}$ precisely parallels the action of $\pi$ on $T^{\uparrow\beta}$.

If this is accomplished, then the action of $\pi^*$ fixes $S$ and so fixes $$\{\{g\}_\beta:g \in G\}_\alpha,$$ while at the same
time $(\pi^*_\beta)_\gamma$ agrees with $\pi_\gamma$ at $g$.  This implies that
$$ \pi_\gamma(g) = (\pi^*_\beta)_\gamma(g) \in (\pi^*_\beta)_\gamma `` G = G. $$
The same construction can be done for $\pi^{-1}$ in place of $\pi$ (for each $g$),
so $\pi$ fixes $G_\beta$.

We construct the allowable permutation $\pi^*$ by the Freedom of Action Theorem (\ref{thm:foa_actions}), by defining the following $\alpha$-action $\xi$.

For every $(\{x\},A) \in {\tt rng}(S)$, we define $\xi_A(x) = x$, and for every $(N,A) \in {\tt rng}(S)$ for $N$ a near-litter, we define $\xi_A(N) = N$.
Extend $T^{\uparrow\beta}$ to a strong $\beta$-support $T^*$ {(as in proposition \ref{prop:extend_to_strong_support})}, then if $(\{x\},B) \in {\tt rng}(T^*)$, we define $(\xi_\beta)_B(x) = \pi_B(x)$, and if $(N,B) \in {\tt rng}(T^*)$ for $N$ a near-litter, we similarly define $(\xi_\beta)_B(N) = \pi_B(N)$.

It now suffices to check that this action is coherent; then, by the Freedom of Action Theorem, we will obtain $\pi^*$ with the property that $\pi^*[S] = S$ and $\pi^*_\beta[T^*] = \pi[T^*]$, so in particular, $(\pi^*_\beta)_\gamma[T] = \pi_\gamma[T]$.

If $\xi_A(N)$ is defined, we have two cases: either $(N,A)$ appears in $S$, in which case $\xi_A(N) = N$, or $B = A\setminus\{\alpha\}$ has maximum element $\beta$ and $(N,B)$ occurs in $T^*$ with $\xi_A(N) = \pi_B(N)$.
Note that if both of these cases are true at the same time, then as $\pi$ fixes $S_{(\beta)}$, the two defined images of $N$ coincide.

In the first case, the case for flexibility is trivial, and if $N \sim f_{\delta,{\tt min}(A_1)}(x,U)$ with $\delta < {\tt min}(A_2)$, then as $S$ is a strong support, the range of $U^{\uparrow A_2}$ is a subset of the range of $S$, as required.

We now consider the second case.
If $N$ is $A$-flexible, we claim that $\pi_B(N)$ is also $A$-flexible.
Suppose that $\pi_B(N) \sim f_{\delta,\gamma}(x, U)$ for $\gamma = {\tt min}(A_1)$ and $\delta < {\tt min}(A_2)$.
Then by the coherence condition, assuming $B$ has length at least 3,
$$ N \sim \pi_B^{-1}``f_{\delta,\gamma}(x, U) \sim f_{\delta,\gamma}((\pi_{B_2}^{-1})_\delta(x), (\pi_{B_2}^{-1})_\delta[U]), $$
so $N$ is $A$-inflexible, giving a contradiction.
In the case that $B$ does not have length at least 3, it must have length exactly 2 because $A$ has length at least 3, so we have $A = \{ \alpha, \beta, -1 \}$ and $B = \{ \beta, -1 \}$.
Note that $(N, B)$ occurs in $T^*$, so it either came from a support condition in $T$, which is impossible as $B = \{ \beta, -1 \}$ does not contain $\gamma$, or it came from a type-raised copy of an $\eta$-support we needed to include in $T^*$ to make it a strong support.
But $\eta \neq -1$ as $-1$-supports cannot contain near-litter support conditions, and if $\eta$ were a proper type index, we would have $-1 < \eta < \beta$ and $\eta \in B$, also giving a contradiction.

Finally, suppose that $N$ is $A$-inflexible, so $N \sim f_{\delta,\gamma}(x, U)$ for $\gamma = {\tt min}(A_1)$ and $\delta < {\tt min}(A_2)$.
By the argument of the previous paragraph, we must conclude that $A$ is of length at least 4, and so $\gamma = {\tt min}(B_1)$ and $\delta < {\tt min}(B_2)$.
As $T^*$ is strong, it contains a suitably type-raised copy of $U$, and so we are done by definition of the action of $\xi$ on $T^*$.

\end{proof}

This completes the proof.  In the formal proof in Lean, what is actually done is a proof that each of the assertions in the finite axiomatization of Hailperin in the version discussed in subsection \ref{ss:hailperin}  holds in all typed versions in our structure for the language of TTT, so it is in fact a model of TTT.  The axioms in Hailperin other than the axiom of type lowering are predicative comprehension axioms and
admit demonstration by the methods of section \ref{ss:predicative_ttt}, done explicitly without metamathematics.  In the formalization, the axiom of type lowering, which contains rather more content than the axiom of set union restricted to sets of singletons which is proved here, is proved by first proving the existence of an iterated image under elementwise application of the singleton operation of the desired set, whose definition is predicative, then repeatedly applying the result of this section that sets of singletons have unions.

\newpage

\section{Conclusions, extended results, and questions}

In a series of subsections, we discuss further results supported by our model construction and further questions.

\subsection{Ways in which this resolution of the NF consistency problem is unexciting}

This is in some ways a rather boring resolution of the NF consistency problem.

NF has no locally interesting combinatorial consequences.   Any stratified fact about sets of a bounded standard size which holds in ZFC will continue to hold in models constructed using this strategy with the parameter $\kappa$ chosen large enough.
That the continuum can be well-ordered or that the axiom of dependent choices can hold, for example, can readily be arranged.  Any theorem about familiar objects such as real numbers which holds in ZFC can be relied upon to hold in our models
(even if it requires Choice to prove), and any situation which is possible for familiar objects is possible in models of NF:  for example, the Continuum Hypothesis can be true or false.  It cannot be expected that NF proves any strictly local stratified result about familiar mathematical objects which is not also a theorem of ZFC.

\subsection{Questions and a result about global choice-like statements}

Questions of consistency with NF of global choice-like statements such as ``the universe is linearly ordered"  or the prime ideal theorem cannot be resolved by the method used here (at least, not without major changes).  Our models falsify the existence of a global linear order at the outset:  a litter cannot be linearly ordered.  We see no reason to believe that the universe cannot be linearly ordered in NF or that the prime ideal theorem cannot hold, but we do not know how to model these situations.

We do have one quite interesting theorem which appears to show that most practical applications of the axiom of choice are supported in our models of TTT and so with care in the resulting models of NF.

Note that in our models of TTT, any relation which the internal theory thinks is a well-ordering actually is a well-ordering.
The reason is that in the metatheory any linear order which is not a well-ordering admits a countable set with no minimal element (the range of a strictly decreasing sequence).  This requires choice, but choice holds in the metatheory.  So, if a relation represented in the model of TTT is not a well-ordering from the standpoint of the metatheory, there is a countable subcollection of its domain which has no minimal element, and this countable subset is actually realized in the model of TTT, because our models of TTT are countably complete, so the model of TTT sees that the relation is not a well-ordering.

Essentially the same argument shows that any relation which a model of TTT constructed by our method thinks is well-founded is actually well-founded in the sense of the metatheory.

\begin{description}

\item[Theorem:]  In our models of TTT, the power set of any well-orderable set can be well-ordered.

\item[Proof:]  Suppose that $X$ is a set in a model of TTT constructed by our method which admits a well-ordering (defined along a particular type sequence).
The well-ordering $\leq_X$ of $X$ which we are given must have a support.  A permutation which fixes the support of $\leq_X$ [in the appropriate sense defined in the paper] will fix $\leq_X$ and, because $\leq_X$ is a well-ordering in the sense of the metatheory, its action on the type of the elements of $X$ (in the type sequence with respect to which $\leq_X$ is a well-ordering) must fix each element of $X$,
so in fact we have a common support for all elements of $X$, which in turn gives a support for each subcollection of $X$ in the sense of the metatheory as implemented in the appropriate type in the particular type sequence, so all of these are sets in the model of TTT (in the appropriate type in the particular type sequence) and further gives a support for the well-ordering of the power set of $X$ induced by any well-ordering of the power set of $X$ in the sense of the metatheory, so the power set of $X$ is well-ordered in the model of TTT.

It is worth remarking that this theorem has been proved in the Lean formalization.

\end{description}

It is an interesting incidental remark that the theorem is equivalent in NF to the assertion that a single big power set, the power set of the set of ordinals, can be well-ordered, and this is seen to hold in models of NF constructed from our models of TTT.

This form of choice seems to support all practical applications of choice in mathematics outside of set theory, and a lot of the applications inside set theory.

Applying the same technique to sets with well-founded extensional relations on them can be used to show that all $\beth$ numbers are alephs, and that the internal representation of an initial segment of the cumulative hierarchy using equivalence classes of well-founded extensional relations with top which is possible in TTT will agree precisely with the cumulative hierarchy in the metatheory as far as it goes (the internal representation will always give a proper initial segment of the hierarchy in the metatheory, of course).

\subsection{Strong unstratified axioms for NF derived from strong partition properties}

NF with strong axioms such as the Axiom of Counting (introduced by Rosser in \cite{rosser}, an admirable textbook based on NF), the Axiom of Cantorian Sets (introduced in \cite{henson})\footnote{Getting Cantorian Sets or Large Ordinals to hold is sensitive to the relationship between $\kappa$ and $\lambda$.  It requires the hypothesis that $\kappa<\lambda=\mu$ to avoid outright refutability of this axiom in resulting models of NF, and then a large cardinal hypothesis.  We note that for the moment we need to use Henson's original formulation of the Axiom of Cantorian Sets, that any {\em well-orderable\/} cantorian set is strongly cantorian;  the situation for non-well-orderable sets is not yet clear to us.} or my\footnote{The first author is speaking.} axioms of Small Ordinals and Large Ordinals (introduced in  my \cite{mybook}, which pretends to be a set theory textbook based on NFU) can be obtained by choosing $\lambda$ large enough to have strong partition properties, more or less exactly as I report in my paper \cite{strongaxioms} on strong axioms of infinity in NFU:  the results in that paper are not all mine, and I owe a good deal to Robert Solovay in that connection (unpublished conversations and \cite{nfub}).

That NF has $\alpha$-models for each standard ordinal $\alpha$ follows by the same methods Jensen used for NFU in his original paper \cite{nfu}:  the proof we give in subsection \ref{subsection:alpha_models} for the existence of $\alpha$-models is essentially the same as Jensen's, and a model of TTT constructed by our method will implement each ordinal $\alpha$ of the metatheory which is less than $\kappa$ (and we can choose $\kappa$ larger than any fixed $\alpha$ and $\lambda$ with cofinality greater than $|\alpha|^{|\alpha|}$).

  No model of NF can contain all countable subsets of its domain;  all well-typed combinatorial consequences
of closure of a model of TST under taking subsets of size $<\kappa$ will hold in our models, but the application of compactness which gets us from TST + Ambiguity to NF forces the existence of externally countable proper classes, a result which has long been known and which also holds in NFU.

\subsection{Some esoteric questions answered}

We mention some esoteric problems which our approach solves.  The Theory of Negative Types of Hao Wang (TST with all integers as types, proposed in \cite{tnt})  has $\omega$-models;  an $\omega$-model of NF gives an $\omega$-model of the theory of negative types  immediately.  The question of existence of $\omega$-models of the theory of negative types was open.

In ordinary set theory, the Specker tree of a cardinal is the tree in which the top is the given cardinal, the children of the top node  are the preimages of the top under the map $(\kappa \mapsto 2^{\kappa})$, and the part of the tree
below each child is the Specker tree of the child.  Forster proved using a result of Sierpi\'nski that the Specker tree of a cardinal must be well-founded (a result which applies in ordinary set theory or in NF(U), with some finesse in the definition of the exponential map in NF(U)).  Given Choice, there is a finite bound on the lengths of the branches in any given Specker tree.  Of course by the Sierpi\'nski result a Specker tree can be assigned an ordinal rank.  The question which was open
was whether existence of a Specker tree of infinite rank is consistent.  It is known that in NF with the Axiom of Counting the Specker tree of the cardinality of the universe is of infinite rank.   Our methods show the existence of a model of NF with the Axiom of Counting.  This gives us a model of TST in which the Specker tree of some cardinal is of infinite rank.  It is straightforward to produce a model of bounded Zermelo set theory either with atoms or with extensionality and failure of foundation from a model of TST, which will satisfy the same stratified sentences:  identify the elements of type 0 either as atoms or with their own singletons, then implement elements of each type $i+1$ as subsets of type $i$ as guided by the original model of TST, identifying type $i+1$ objects with elements of lower types which have been assigned the same extension.  This shows that we can get a model of bounded Zermelo set theory with atoms or without foundation in which there is a cardinal with Specker tree of infinite rank.  It is straightforward with our methods to get a model of TTT with the Axiom of Counting with inaccessibles cofinal in $\lambda$, and in this model we will have an inaccessible greater than a cardinal with Specker tree of infinite rank,
and so a model of ZFA in which there is a cardinal with Specker tree of infinite rank.  It's actually a bit easier than this:  the existence of a cardinal with Specker tree of infinite rank holds in any model built by our method in which $\lambda>\omega$, but referring back to the result about NF + Counting makes it easier to give a brief account here.

The further question remains as to whether it is consistent with ZF that there is a cardinal with Specker tree of infinite rank.  We would guess that this could be shown to be true using symmetric forcing extensions, but we are not expert at this and leave the question to others.

\subsection{The question of consistency strength}

We believe that NF is no stronger than TST + Infinity, which is of the same strength as Zermelo set theory with separation restricted to bounded formulas (\cite{kemeny}).  Our work here does not show this, as we need enough Replacement for
existence of $\beth_{\omega_1}$ at least.  We leave it as an interesting further task, possibly for others, to tighten things up and show the minimal strength that we expect holds.

\subsection{The wide open question:  what do models of NF look like in general?}

Another question of a very general and amorphous nature which remains is:  what do models of NF (or TTT) look like in general?  Are all models of NF in some way like the ones we describe, or are there models of quite a different character?  There are very special assumptions which we made by fiat in building our model of TTT which do  not seem at all inevitable in general models of this theory.

\subsection{Postscript}

Inevitably, philosophical issues come up in connection with a system of set theory proposed by a philosopher.  I get a lot of congratulations for vindicating Quine's foundational agenda, but in fact this is not part of my purpose here.  It is not even clear to me that Quine {\em had\/} a foundational agenda in which his technical proposal of this set theory had a special place.

The results of this paper show that if one really wanted to, one could use NF as a foundation for mathematics.  It is odd that it disproves Choice, but the principle ``power sets of well-orderable sets are well-orderable", which holds in models constructed in this way, supports most applications of Choice.

Our opinion is that Quine's proposal of NF was based on a mistake.  He discusses whether to assume strong extensionality in the original paper, and his explanation of reasons for choosing strong extensionality contains an actual mathematical error.
We believe that the correct system to propose was NFU, and if he had proposed NFU, the history of this kind of set theory might have been different.

NFU is serviceable as a foundation for mathematics, and consistent with Choice and various strong axioms of infinity.  It is odd that NFU + Choice proves that there are urelements (the same odd fact that NF disproves Choice), but no more than odd.
The type discipline of stratification is something that one must work to get used to, and it has been remarked that working with indexed families of sets is extremely awkward in NFU (and would be similarly awkward in NF).  The view of the world which NFU supports is basically the same as that of ZFC:  the natural models of NFU are obtained by considering an initial segment of the cumulative hierarchy with an external automorphism which moves a level, and NFU can interpret discussion of exactly such a structure internally.

Unlike NF, NFU on introspection can tell one quite a lot about what its models should look like (as ZFC can with its own awareness of initial segments of the cumulative hierarchy; in NFU, analysis of the isomorphism classes of well-founded extensional relations gives both an interpretation of an initial segment of the Zermelo style universe and an interpretation of NFU itself if one has strong enough assumptions).    NF tells one very little about what its intended world is like (it can, being an extension of NFU, internally construct a lot of information about an interpretation of NFU with lots of urelements, but there is not an obvious way to find out how an extensional world is constructed from internal evidence in NF itself).

So, we do not believe this paper is a philosophical milestone.  If there was one, it happened in 1969 when Jensen showed that NFU is consistent, and nobody noticed.

We do believe that there are interesting questions to investigate about NF.  Our paper does not settle all the questions about this system which have developed in the minds of people who have worked with it since 1937.  In fact, the construction is based on special assumptions  and does not seem to give much of an idea of what general models of this theory might look like.

It might be viewed as philosophically interesting that this proof was formally verified.  That is really an advertisement for a quite different foundational system, the logic of the Lean proof verification system.   That it needed to be formally verified I believe reflects an interesting complexity in the mathematics here, not only the deficiencies of the first author as an expositor.

\end{document}